\documentclass[11pt]{article}
\usepackage[english]{babel}
\usepackage[a4paper,scale=0.75,hcentering]{geometry}
\usepackage[hyperfootnotes=false]{hyperref} 
\usepackage{setspace}
\hypersetup{
	colorlinks = true,
	linkcolor = {blue},
	urlcolor = {red},
	citecolor = {blue}
}
\usepackage{amsmath, amsthm, amssymb, mathrsfs, mathtools, enumitem}
\usepackage[x11names]{xcolor}
\allowdisplaybreaks

\newcounter{results}[section] 

\theoremstyle{plain}
\newtheorem{theorem}[results]{Theorem}
\newtheorem{lemma}[results]{Lemma}
\newtheorem{proposition}[results]{Proposition}
\newtheorem{corollary}[results]{Corollary}

\newtheorem*{theorem*}{Theorem}
\newtheorem*{lemma*}{Lemma}
\newtheorem*{proposition*}{Proposition}
\newtheorem*{corollary*}{Corollary}
\newtheorem*{exercise*}{Exercise}
\newtheorem*{fact*}{Fact}
\newtheorem*{problem*}{Problem}
\newtheorem*{conjecture*}{Conjecture}

\theoremstyle{remark}
\newtheorem{remark}[results]{Remark}

\newtheorem*{remark*}{Remark}
\newtheorem*{question*}{Question}

\theoremstyle{definition}

\newtheorem*{definition*}{Definition}
\newtheorem*{example*}{Example}

\numberwithin{equation}{section}

\newcommand{\N}{\ensuremath{\mathbb N}} 
\newcommand{\R}{\ensuremath{\mathbb R}} 

\renewcommand{\S}{\ensuremath{\mathbb S}} 
\renewcommand{\P}{\ensuremath{\mathcal P}}

\newcommand{\M}{\ensuremath{\mathcal M}}
\makeatletter 
\DeclarePairedDelimiter{\@tmpabs}{\lvert}{\rvert}
\newcommand{\@absstar}[1]{{\@tmpabs*{#1}}}
\newcommand{\@absnostar}[2][]{{\@tmpabs[#1]{#2}}}
\newcommand{\abs}{\@ifstar\@absstar\@absnostar}
\makeatother

\makeatletter 
\DeclarePairedDelimiter{\@tmpnorm}{\lVert}{\rVert}
\newcommand{\@normstar}[1]{{\@tmpnorm*{#1}}}
\newcommand{\@normnostar}[2][]{{\@tmpnorm[#1]{#2}}}
\newcommand{\norm}{\@ifstar\@normstar\@normnostar}
\makeatother

\newcommand{\weakto}{\ensuremath{\rightharpoonup}} 
\newcommand{\weakstarto}{\ensuremath{\stackrel{*}{\rightharpoonup}}} 

\newcommand{\de}{\ensuremath{\,\mathrm d}}

\renewcommand{\P}{\ensuremath{\mathcal{P}}}
\newcommand{\PP}{\ensuremath{\mathring{\mathcal{P}}}}

\DeclareMathOperator{\supp}{supp}

\ifdefined\comma
	\renewcommand{\comma}{\ensuremath{\, \text{, }}}
\else
	\newcommand{\comma}{\ensuremath{\, \text{, }}}
\fi

\usepackage[nameinlink,capitalise,sort]{cleveref} 
\crefname{equation}{}{} 
\crefname{enumi}{}{} 

\begin{document}

\title{\bfseries Improved Sobolev and Moser--Trudinger--Onofri Inequalities under Constraints:
Sharpness and Quantitative Remainders\thanks{This work is supported by the National Key R\&D Program of China (Grant 2023YFA1010001) and NSFC (12571123). E-mail addresses: zhou-yx22@mails.tsinghua.edu.cn (Zhou), zou-wm@mail.tsinghua.edu.cn (Zou).}}

\author{{\bfseries Yuxuan Zhou, Wenming Zou}\\ {\footnotesize \itshape Department of Mathematical Sciences, Tsinghua University, Beijing 100084, China.}}

\date{}

\maketitle

\begin{abstract}
{\small
In this paper, we study improved Sobolev and
Moser--Trudinger--Onofri inequalities on the sphere under moment
constraints, with particular emphasis on sharp leading constants and
quantitative remainder terms.

For the Sobolev inequality, we sharpen the almost sharp inequality of
Hang and Wang \cite{Hangwang} (Int. Math. Res. Not., 2022) to a sharp
one at $p=2$ under vanishing moment constraints of orders $m=1,2,3$.
We also prove a general failure principle: if the unconstrained sharp
Sobolev inequality fails, then the corresponding sharp improved
inequality under moment constraints fails as well.

For the Moser--Trudinger--Onofri inequality, we sharpen the almost sharp
inequalities of Chang and Hang \cite{Cha2} and Hang \cite{Hang}
(Comm. Pure Appl. Math., 2022) to a sharp improved inequality at the
critical order $s=\frac{n}{2}$ under arbitrary finite-order moment
constraints. We further investigate the best remainder constant. In
contrast with the two-dimensional first-moment case of Gui and
Moradifam \cite{Guimo} (Invent. Math., 2018), where the best remainder
vanishes, we show that it is positive in several higher-dimensional
regimes.

We also obtain quantitative remainder estimates for general almost
sharp improved Sobolev and Moser--Trudinger--Onofri inequalities, which
in turn yield sharp inequalities at different scales. We further extend
recent second-moment refinements to higher dimensions and higher orders.

Finally, we obtain analogous results on sharpness, non-sharpness, and
quantitative remainder estimates for improved Sobolev-type inequalities
under symmetry constraints on general smooth closed manifolds.



}
\end{abstract}
\tableofcontents

\vskip0.6in

\section{Introduction}
\subsection{Background and motivation}
Assume that $n\in\N_+$, $1\leq p<n$, and set $p^*=\frac{np}{n-p}$.
The classical Sobolev inequality states that for any function
$u\in W^{1,p}(\R^n)$,
\begin{align}\label{Sob0}
    \norm{u}_{L^{p^*}(\R^n)} \leq S_{n,p}\norm{\nabla u}_{L^p(\R^n)},
\end{align}
where $S_{n,p}$ denotes the optimal Sobolev constant. It follows from
\cite{Aub1,Tal} that
\begin{align}
    S_{n,1}=\frac{1}{n}\left(\frac{n}{\omega_{n-1}}\right)^\frac{1}{n},\nonumber
\end{align}
and for $p>1$,
\begin{align}
    S_{n,p}
    = \frac{1}{n}\left(\frac{n(p-1)}{n-p}\right)^{1-\frac{1}{p}}
    \left(\frac{n!}{\Gamma\left(\frac{n}{p}\right)
    \Gamma\left(n+1-\frac{n}{p}\right)\omega_{n-1}}\right)^{\frac{1}{n}},\nonumber
\end{align}
where $\omega_{n-1}$ denotes the volume of $\S^{n-1}$ and
$\Gamma(s):=\int_0^{\infty} t^{s-1}e^{-t}\,\mathrm{d}t$ for $s>0$.

The Moser--Trudinger inequality, which can be viewed as the endpoint case
of \eqref{Sob0}, states that, for any bounded open set
$\Omega\subset \R^n$ with $n\geq 2$ and any function
$u\in W^{1,n}_0(\Omega)\setminus\{0\}$,
\begin{align}\label{MT0}
    \int_\Omega
    \exp\left(
    a_n\frac{|u|^{\frac{n}{n-1}}}
    {\norm{\nabla u}_{L^n(\Omega)}^{\frac{n}{n-1}}}
    \right)\,\mathrm{d}x
    \leq C(n)|\Omega|.
\end{align}
Here $a_n = n\omega_{n-1}^{\frac{1}{n-1}}$ is the optimal constant; see \cite{Moser,Tru}. In \cite{Ada}, Adams generalized \eqref{MT0} to higher-order Sobolev spaces.
Assume that $s,n\in\N_+$ and $1\leq s<n$. For any bounded open set
$\Omega\subset \R^n$ and any function
$u\in W^{s,\frac{n}{s}}_0(\Omega)\setminus\{0\}$,
\begin{align}\label{MT1}
    \int_\Omega
    \exp\left(
    a_{s,n}\frac{|u|^{\frac{n}{n-s}}}
    {\norm{\Delta^{\frac{s}{2}} u}_{L^{\frac{n}{s}}(\Omega)}^{\frac{n}{n-s}}}
    \right)\,\mathrm{d}x
    \leq C(n,s)|\Omega|.
\end{align}
The optimal constant $a_{s,n}$ is given by
\begin{align}
    a_{s,n}
    = \frac{n}{\omega_{n-1}}
    \left(
    \frac{\pi^{\frac{n}{2}}2^s
    \Gamma\left(\lfloor\frac{s+1}{2}\rfloor\right)}
    {\Gamma\left(\frac{n}{2}-\lfloor\frac{s}{2}\rfloor\right)}
    \right)^{\frac{n}{n-s}}.\nonumber
\end{align}
When $s$ is odd, $\Delta^{\frac{s}{2}}u$ should be understood as
$\nabla\Delta^{\frac{s-1}{2}}u$.

\vskip0.12in

In this paper, we are concerned with generalizations of
\eqref{Sob0}, \eqref{MT0}, and \eqref{MT1} to general smooth, closed, connected manifolds.
Let $(M,g)$ be a smooth closed manifold of dimension $n\geq 2$.
For $1\leq p<n$, the Sobolev inequality for functions
$u\in W^{1,p}(M)$ takes the form
\begin{align}\label{AB}
    \norm{u}_{L^{p^*}(M)}
    \leq A\norm{\nabla u}_{L^p(M)} + B\norm{u}_{L^{p}(M)}.
\end{align}
Unlike the Euclidean case \eqref{Sob0}, it is easy to see that the term
$\norm{u}_{L^{p}(M)}$ on the right-hand side of \eqref{AB} cannot be removed.
The problem of optimizing the coefficients $A$ and $B$ in \eqref{AB}
is known as the \textbf{$AB$-program}, and a detailed introduction
can be found in the monographs \cite{Aub3,Heb0,Heb1}.

\vskip0.12in

In the pioneering work \cite{Aub1}, Aubin showed that for any
$\varepsilon>0$ there exists a positive constant
$B(M,p,\varepsilon)>0$ such that
\begin{align}\label{Sob2}
    \norm{u}_{L^{p^*}(M)}^p
    \leq (S_{n,p}^p+\varepsilon)\norm{\nabla u}_{L^p(M)}^p
    + B(M,p,\varepsilon)\norm{u}_{L^{p}(M)}^p.
\end{align}
Moreover, this inequality is \textbf{almost sharp} in the sense that
the coefficient $S_{n,p}^p+\varepsilon$ cannot be replaced by any
number smaller than $S_{n,p}^p$.
It is a delicate issue whether one can take $\varepsilon=0$ in
\eqref{Sob2} and obtain the following \textbf{sharp} inequality:
\begin{align}\label{Sob3}
     \norm{u}_{L^{p^*}(M)}^p
     \leq S_{n,p}^p\norm{\nabla u}_{L^p(M)}^p
     + B(M,p)\norm{u}_{L^{p}(M)}^p.
\end{align}
On $\S^{n}$, it was shown in \cite{Aub1}
(see also \cite{Heb2}) that \eqref{Sob3} holds for every $1\leq p\leq 2$.
However, it follows from \cite{Dru} that \eqref{Sob3} cannot hold
for $2<p<\sqrt{n}$.
We refer interested readers to \cite[Section~4]{Heb1}
for further related results.

\vskip0.12in

For the higher-order Moser--Trudinger inequality, it follows from
\cite{Bra,Fon} that the following \textbf{sharp} inequality
\begin{align}\label{MT3}
    \int_M
    \exp\left(
    a_{s,n}\frac{|u|^{\frac{n}{n-s}}}
    {\norm{\Delta^{\frac{s}{2}} u}_{L^{\frac{n}{s}}(M)}^{\frac{n}{n-s}}}
    \right)\,\mathrm{d}\mu
    \leq c(M,s)
\end{align}
holds for any
$u\in W^{s,\frac{n}{s}}(M)\setminus\{0\}$
satisfying $\int_M u\,\mathrm{d}\mu=0$. Here $\mu$ is the measure associated with the metric $g$. Combining \eqref{MT3} with Young's inequality, we obtain
the higher-order Moser--Trudinger--Onofri inequality
(see \cite{Hang,Moser,Onofri})
\begin{align}\label{MTO0}
    \log\int_M e^{nu}\,\mathrm{d}\mu
    \leq \alpha_{s,n}
    \norm{\Delta^{\frac{s}{2}} u}_{L^{\frac{n}{s}}(M)}^{\frac{n}{s}}
    + C(M,s).
\end{align}
Here
\begin{align}
    \alpha_{s,n}
    = s\left(\frac{(n-s)\omega_{n-1}}{n}\right)^{\frac{n-s}{s}}
    \left(
    \frac{\Gamma\left(\frac{n}{2}-\lfloor\frac{s}{2}\rfloor\right)}
    {\pi^{\frac{n}{2}}2^s
    \Gamma\left(\lfloor\frac{s+1}{2}\rfloor\right)}
    \right)^{\frac{n}{s}}\nonumber
\end{align}
is the optimal constant.
A careful analysis of the constants appearing in
\eqref{MT3} and \eqref{MTO0} can be found in \cite{Yang}.

\vskip0.12in

In the seminal work \cite{Aub2}, Aubin showed that by imposing
vanishing first-order moment constraints, one can improve
the constants in both the Sobolev inequality and the
first-order Moser--Trudinger--Onofri inequality on $\S^{n}$.
More precisely, for any $\varepsilon>0$ and any function
$u\in W^{1,p}(\S^n)$ satisfying
\begin{align}
    \int_{\S^n} x_i |u|^{p^*}\,\mathrm{d}\mu = 0,
    \qquad 1\leq i\leq n+1,
    \nonumber
\end{align}
there exists a positive constant $B(n,p,\varepsilon)$ such that
\begin{align}\label{Sob4}
    \norm{u}_{L^{p^*}(\S^n)}^p
    \leq \left(2^{-\frac{p}{n}} S_{n,p}^p + \varepsilon\right)
    \norm{\nabla u}_{L^p(\S^n)}^p
    + B(n,p,\varepsilon)\norm{u}_{L^{p}(\S^n)}^p .
\end{align}
Similarly, for any $\varepsilon>0$ and any function
$u\in W^{1,n}(\S^n)$ satisfying
\begin{align}
    \int_{\S^n} x_i e^{n u}\,\mathrm{d}\mu = 0,
    \qquad 1\leq i\leq n+1,
    \nonumber
\end{align}
there exists a positive constant $C(n,\varepsilon)$ such that
\begin{align}\label{MTO2}
    \log\int_{\S^n} e^{n u}\,\mathrm{d}\mu
    \leq \left(\frac{\alpha_{1,n}}{2} + \varepsilon\right)
    \norm{\nabla u}_{L^{n}(\S^n)}^n
    + C(n,\varepsilon).
\end{align}
Moreover, the inequalities \eqref{Sob4} and \eqref{MTO2}
are \textbf{almost sharp} in the same sense as \eqref{Sob2}.
These inequalities play an important role in the analysis
of prescribing scalar curvature problems on $\S^n$;
see \cite{Aub2,Aub3,Cha3,Cha4,Heb1} and the references therein
for further discussions.

\vskip0.1in

Recently, generalizations of \eqref{Sob4} and \eqref{MTO2}
incorporating higher-order moment constraints and
higher-order Moser--Trudinger--Onofri inequalities
have been established in a series of works
\cite{Cha2,Hang,Hangwang}.
To state these results, we first introduce some notation.

\vskip0.1in

For $m,n\in\N_+$, we set
\begin{align}
    \P_{m,n}
    =\{\text{all polynomials on $\R^{n+1}$ of degree at most $m$}\}
\end{align}
and
\begin{align}
    \PP_{m,n}
    =\left\{f\in \P_{m,n}:\;
    \int_{\S^n} f\,\mathrm{d}\mu = 0 \right\}.
\end{align}
Let $\Delta(\S^n)$ denote the space of probability measures
on $\S^n$ supported on countably many points.
For $\nu=\sum\limits_i \nu_i \delta_{\xi_i}\in\Delta(\S^n)$
and $0\leq \theta<1$, we define
$\norm{\nu}_\theta := \sum\limits_i \nu_i^\theta$.

The following two notions, which originate from numerical
analysis and scientific computing
(see, for instance, \cite{de,Del,Put}),
are useful for formulating improved inequalities
under vanishing higher-order moment constraints:
\begin{align}
    \M^c_{m,n}
    =\left\{\nu\in\Delta(\S^n):\;
    \int_{\S^n} f\,\mathrm{d}\nu = 0
    \;\text{for all } f\in\PP_{m,n}\right\},
\end{align}
\begin{align}\label{Cub1}
    \Theta(m,\theta,n)
    =\inf\left\{\norm{\nu}_\theta:\;
    \nu\in \M^c_{m,n}\right\}.
\end{align}
The elements of $\M^c_{m,n}$ are usually referred to as cubature formulas of strength $m$. Estimates for $\Theta(m,\theta,n)$, particularly in
the case $\theta=0$, have been widely studied. It is known that
$\Theta(1,\theta,n)=2^{1-\theta}$, $\Theta(2,\theta,n)=(n+2)^{1-\theta}$, and $\Theta(3,\theta,n)=(2n+2)^{1-\theta}$. We refer interested readers to
\cite{de,Del,Hangwang,Put} and the references therein
for a more detailed discussion of cubature formulas,
spherical designs, and estimates for $\Theta$.
\vskip0.12in

We now state the \textbf{almost sharp} improved Sobolev
and Moser--Trudinger--Onofri inequalities
under vanishing higher-order moment constraints.
In \cite[Theorem~1.1]{Hangwang},
Hang and Wang showed that for $m,n\in\N_+$, $1<p<n$,
$\varepsilon>0$, and any function
$u\in W^{1,p}(\S^n)$ satisfying
\begin{align}\label{or1}
    \int_{\S^n} f |u|^{p^*}\,\mathrm{d}\mu = 0
    \quad \text{for all } f\in \PP_{m,n},
\end{align}
the following \textbf{almost sharp} improved
Sobolev inequality holds:
\begin{align}\label{main1}
    \norm*{u}_{L^{p^*}(\S^n)}^p
    \leq
    \left(
        \frac{S_{n,p}^p}
        {\Theta\left(m,\frac{n-p}{n},n\right)}
        +\varepsilon
    \right)
    \norm*{\nabla u}_{L^p(\S^n)}^p
    + B(m,n,p,\varepsilon)
    \norm*{u}_{L^{p}(\S^n)}^p,
\end{align}
for some positive constant $B(m,n,p,\varepsilon)$.
\vskip0.12in

For the higher-order Moser--Trudinger--Onofri inequality,
it follows from Chang--Hang
\cite[Theorem~1.2]{Cha2}
and Hang
\cite[Theorems~3.2, 4.5, 4.11]{Hang}
that for $s,m,n\in\N_+$, $1\leq s<n$,
$\varepsilon>0$, and any function
$u\in W^{s,\frac{n}{s}}(\S^n)$ satisfying
$\int_{\S^n} u\,\mathrm{d}\mu = 0$ and
\begin{align}\label{or2}
    \int_{\S^n} f e^{n u}\,\mathrm{d}\mu = 0
    \quad \text{for all } f\in \PP_{m,n},
\end{align}
the following \textbf{almost sharp} improved inequality holds:
\begin{align}\label{main2}
    \log \int_{\S^n} e^{n u}\,\mathrm{d}\mu
    \leq
    \left(
        \frac{\alpha_{s,n}}{\Theta(m,0,n)}
        +\varepsilon
    \right)
    \norm{\Delta^{\frac{s}{2}} u}_{L^{\frac{n}{s}}(\S^n)}^{\frac{n}{s}}
    + C(m,n,s,\varepsilon),
\end{align}
for some positive constant $C(m,n,s,\varepsilon)$.
Since $\Theta(1,\theta,n)=2^{1-\theta}$
for any $0\leq \theta<1$,
one can readily recover \eqref{Sob4} and \eqref{MTO2}
by taking $m=1$ in \eqref{main1} and \eqref{main2},
respectively.
\vskip0.2in
Motivated by the \textbf{$AB$-program} for the inequality \eqref{AB},
it is natural to ask whether one can take $\varepsilon=0$
in \eqref{main1} and \eqref{main2},
and what can be said about the remainders
$B_m$ and $C_m$.
The study of these questions is not only important
from the perspective of functional inequalities,
but also closely related to various prescribing
curvature problems on $\S^n$.
So far, however, only limited progress has been made.
The only completely solved case is due to
Gui and Moradifam \cite{Guimo}.
In \cite{Guimo}, they applied the sphere covering inequality
to give a complete analysis of \eqref{main2}
in the case $m=s=1$ and $n=2$,
thereby confirming a conjecture proposed in \cite{Cha3}.
More recently, Gui, Li, Wei, and Ye \cite{Guili1} generalized the results in \cite{Guimo} to Beckner's inequality involving the GJMS operator (see also \cite{Guihu1,Guihu2,Guili,Guili2,Li} for partial results). Since the GJMS energy is strictly larger than the fractional Sobolev energy when $n\geq 3$, their results do not imply that one can take $\varepsilon=0$ in \eqref{main2} when $s=\frac{n}{2}$.

\vskip0.12in

It is also worth mentioning that the recent proofs of \eqref{main1} and \eqref{main2} in
\cite{Cha2,Hang,Hangwang} are based on careful use of concentration--compactness principles
(see, for instance,
\cite{Hang,Lions1,Lions2}). These compactness-based arguments
not only allow one to handle vanishing higher-order
moment constraints and yield concise and transparent proofs,
but are also highly flexible and can be adapted to a wide range
of inequalities on general compact manifolds
(see, e.g.,
\cite{An,Chenwei,Chenzhang,Ho,Hu,Sun,Yan1,Yan2}). However, the reliance on concentration--compactness principles
precludes any quantitative information on the remaining constants
$B_m$ and $C_m$ appearing in
\eqref{main1} and \eqref{main2},
let alone the possibility of determining
whether the parameter $\varepsilon$ can be taken to be zero.

\vskip0.12in

In this paper, we make further progress on these problems. By developing a new approach
to \eqref{main1} and \eqref{main2}, we establish sharpness results in several cases. We also provide non-sharp examples through a general failure principle. Moreover, we obtain quantitative estimates for the remainders in general cases, yielding sharp inequalities at different scales. Our approach also applies to improved inequalities under symmetry constraints.

\vskip0.62in

\subsection{Statements of main results}
Let us summarize our main results. Our first result implies that one can take $\varepsilon=0$ in \eqref{main1} when $m\in\{1,2,3\}$, $p=2$, and $n\geq 3$.
\begin{theorem}\label{thm1.1}
    Assume $m\in\{1,2,3\}$, $p=2$, and $n\geq 3$. Then there is a constant $B=B(m,n)<\infty$ such that every $u\in H^1(\S^n)$ satisfying \eqref{or1} also satisfies
    \begin{align}\label{newmain1}
         \norm*{u}_{L^{2^*}(\S^n)}^2
    \leq
        \frac{S_{n,2}^2}
        {\Theta\left(m,\frac{n-2}{n},n\right)}
    \norm*{\nabla u}_{L^2(\S^n)}^2
    + B\norm*{u}_{L^{2}(\S^n)}^2.
    \end{align}
\end{theorem}
To our knowledge, \eqref{newmain1} gives the first \textbf{sharp} improved Sobolev inequality under moment constraints. The restriction $m\in\{1,2,3\}$ arises from a technical issue: our proof relies on properties of extremal measures attaining the value $\Theta$ when $m\in\{1,2,3\}$ (see Remark \ref{rmk2.5} for details).
\vskip0.1in
Our second result provides a general failure principle: if the \textbf{sharp} Sobolev inequality \eqref{Sob3} fails, then the corresponding sharp improved inequality also fails.
\begin{theorem}\label{thm1.2}
Assume $1\leq p<n$ and the \textbf{sharp} Sobolev inequality \eqref{Sob3} fails on $\mathbb S^n$. Then for any $m\in\N_+$, there
does not exist a finite constant $B=B(m,n,p)$ such that
\begin{align}\label{fail}
    \|u\|_{L^{p^*}(\S^n)}^{p}
\leq
\frac{S_{n,p}^{p}}{\Theta\left(m,\frac{n-p}{n},n\right)}
\|\nabla u\|_{L^p(\S^n)}^{p}
+B\|u\|_{L^p(\S^n)}^{p}
\end{align}
holds for every $u\in W^{1,p}(\mathbb S^n)$ satisfying \eqref{or1}.
\end{theorem}
By \cite[Theorem 1]{Dru}, we have the following direct corollary.
\begin{corollary}\label{cor1.3}
    Assume $2<p<\sqrt{n}$. Then for any $m\in\N_+$, there
does not exist a finite constant $B=B(m,n,p)$ such that \eqref{fail}
holds for every $u\in W^{1,p}(\mathbb S^n)$ satisfying \eqref{or1}.
\end{corollary}
To our knowledge, these two results provide the first non-sharp cases of improved Sobolev inequalities under moment constraints.

\vskip0.1in
Our third result implies that one can take $\varepsilon=0$ in \eqref{main2} for $s=\frac{n}{2}$, $n\geq 2$, and every $m\in\N_+$. In particular, we allow $s\notin\N_+$ when $n$ is odd.

\vskip0.22in

\begin{theorem}\label{thm1.4}
    Assume $m,n\in\N_+$, $n\geq 2$, and $s=\frac{n}{2}$. Then there is a constant $C=C(m,n)<\infty$ such that every $u\in H^s(\S^n)$ satisfying $\int_{\S^n} u\,\mathrm{d}\mu = 0$ and \eqref{or2} also satisfies
    \begin{align}\label{newmain2}
        \log\left(\omega_n^{-1} \int_{\S^n} e^{n u}\,\mathrm{d}\mu\right)
    \leq\frac{n}{2(n-1)!\omega_n\Theta(m,0,n)}\norm{(-\Delta)^{\frac{n}{4}} u}_{L^{2}(\S^n)}^{2}
    + C,
    \end{align}
\end{theorem}
To our knowledge, \eqref{newmain2} gives the first \textbf{sharp} improved Moser--Trudinger--Onofri inequalities since the work of Gui and Moradifam \cite{Guimo}, which treats the special case $s=1,n=2,m=1$.
\vskip0.1in
Our next result concerns the best remainder in \eqref{newmain2}, which we denote by $C_{\mathrm{best}}(m,n)$. It follows from \cite{Guimo} that $C_{\mathrm{best}}(1,2)=0$ and that all extremal functions are constant. In recent work \cite{Guili1}, Gui, Li, Wei, and Ye investigated a higher-order Beckner's inequality. They showed that, for the \textbf{sharp} improved Beckner's inequality, the remainder is $0$ and all extremal functions are constant. Based on these results, one might expect that $C_{\mathrm{best}}(m,n)=0$ for general $m,n$. However, we show that this is false, at least for large $n$.
\begin{theorem}\label{thm1.5}
    $C_{\mathrm{best}}(m,n)>0$ in the following four cases: $m=1,n\geq 5$; $m=2,n\geq 12$; $m=3,n\geq 18$; $m\geq 4,n\geq 18(m+1)$.
\end{theorem}
We remark that for $m\geq2$, our lower bound on $n$ is rough. For $m=1$, we devote particular attention to obtaining the lower bound $n\geq 5$, since $C_{\mathrm{best}}(1,2)=0$ is already known. We leave open the interesting question of whether $C_{\mathrm{best}}(1,3)$ and $C_{\mathrm{best}}(1,4)$ vanish.
\vskip0.1in
The following two theorems provide quantitative remainders for \eqref{main1} and \eqref{main2}, respectively, over the full range of parameters.
\begin{theorem}\label{thm1.6}
    Assume $m,n\in\N_+$, $1\leq p<n$, and $0<\varepsilon<1$. Assume further that the \textbf{sharp} Sobolev inequality \eqref{Sob3} holds on
$\S^n$. Then there exists a positive constant $B=B(m,n,p)$ such that, for every $u\in W^{1,p}(\S^n)$ satisfying \eqref{or1},
\begin{align}
    \norm*{u}_{L^{p^*}(\S^n)}^p
    \leq
    \left(
        \frac{S_{n,p}^p}{\Theta(m,\frac{n-p}{n},n)}+\varepsilon
    \right)
    \norm*{\nabla u}_{L^p(\S^n)}^p
    +B\varepsilon^{-\gamma(p)}
    \norm*{u}_{L^p(\S^n)}^p,\nonumber
\end{align}
where
\begin{align}\label{newcst}
    \gamma(p)
    =
    \begin{cases}
        \dfrac{p^2}{2(n-p)}, & 1\leq p\leq2,\\[0.4em]
        \dfrac{np}{2(n-p)}-1, & 2<p<n.
    \end{cases}
\end{align}
\end{theorem}
\begin{theorem}\label{thm1.7}
Assume $s,m,n\in\N_+$, $1\leq s<n$, and $0<\varepsilon<1$. Then there exists a positive constant $C=C(m,n,s)$ such that every function $u\in W^{s,\frac{n}{s}}(\S^n)$ satisfying
$\int_{\S^n} u\,\mathrm{d}\mu = 0$ and \eqref{or2} also satisfies
\begin{align}
     \log \int_{\S^n} e^{n u}\,\mathrm{d}\mu
    \leq
    \left(
        \frac{\alpha_{s,n}}{\Theta(m,0,n)}
        +\varepsilon
    \right)
    \norm{\Delta^{\frac{s}{2}} u}_{L^{\frac{n}{s}}(\S^n)}^{\frac{n}{s}}
    +C\varepsilon^{-\frac{s}{n-s}}.\nonumber
\end{align}
\end{theorem}
We remark that the assumption that \eqref{Sob3} holds in Theorem \ref{thm1.6} is made only for simplicity: when it does not hold, the remainder involves the term $B(m,n,p,\varepsilon)$ in \eqref{main1}. Moreover, we prove stronger results than the two theorems above (see Theorems \ref{thm5.1} and \ref{thm6.2}), in which the orthogonality conditions are replaced by almost-orthogonality conditions.

As a direct consequence of the two theorems above, we obtain \textbf{sharp} improved inequalities at different scales.
\begin{corollary}\label{cor1.8}
    Under the same assumptions as in Theorem \ref{thm1.6}, there exists a positive constant $B=B(m,n,p)$ such that
    \begin{align}\label{scaline3}
        \norm*{u}_{L^{p^*}(\S^n)}^{\Bar{\gamma}(p)}
        \leq
        \Bigg(
            \frac{S_{n,p}}{\Theta\left(m,\frac{n-p}{n},n\right)^{\frac{1}{p}}}
        \Bigg)^{\Bar{\gamma}(p)}
        \norm*{\nabla u}_{L^p(\S^n)}^{\Bar{\gamma}(p)}
        +B\norm*{u}_{L^{p}(\S^n)}^{\Bar{\gamma}(p)},
    \end{align}
    where
    \begin{align}
        \Bar{\gamma}(p)=\frac{p}{1+\gamma(p)}=
        \begin{cases}
            \dfrac{2p(n-p)}{2(n-p)+p^2}, & 1\leq p\leq2,\\[0.4em]
        \dfrac{2(n-p)}{n}, & 2<p<n.
        \end{cases}
    \end{align}
\end{corollary}
\begin{corollary}\label{cor1.9}
    Under the same assumptions as in Theorem \ref{thm1.7}, there exists a positive constant $C=C(m,n,s)$ such that
    \begin{align}\label{scaline2}
        \log^{1-\frac{s}{n}}\Bigg( \frac{1}{\omega_{n}}\int_{\S^n} e^{n u}\,\mathrm{d}\mu\Bigg)
        \leq
        \left(
            \frac{\alpha_{s,n}}{\Theta(m,0,n)}
        \right)^{1-\frac{s}{n}}
        \norm{\Delta^{\frac{s}{2}} u}_{L^{\frac{n}{s}}(\S^n)}^{\frac{n}{s}-1}
        + C.
    \end{align}
\end{corollary}
The inequalities \eqref{scaline3} and \eqref{scaline2} are connected to the \textbf{$AB$-program} for Sobolev-type inequalities at different scales. We refer interested readers to \cite[Section 4.4]{Heb1} for a detailed discussion on this topic.
\vskip0.1in
Finally, we mention two other types of improved inequalities to which our approach can be extended.

\vskip0.1in

The first type comes from the recent work \cite{Chenzhang}. In \cite[Theorem 1.1]{Chenzhang}, Chen and Zhang established
\textbf{almost sharp} improved first-order Moser--Trudinger--Onofri inequalities
under the constraint that the second-order moments deviate suitably from zero
on $\S^2$. This constraint can be viewed as an interesting variant of \eqref{or2}. Here we generalize their results to higher-order inequalities on
$\S^n$ for all $n \geq 2$ and obtain quantitative remainders (see Subsection \ref{sec6.2} for details).

\vskip0.1in

The second type arises from symmetry constraints, which form another classical topic. A brief introduction is given in Section \ref{sec7}. We focus on improved Sobolev and higher-order Moser--Trudinger--Onofri inequalities on general smooth, closed manifolds under symmetry constraints and obtain various new results on sharpness, non-sharpness, and quantitative remainder estimates (see Theorems \ref{thm7.1}, \ref{thm7.2}, and \ref{thm7.8} for details).

\vskip0.36in

\subsection{Outline of strategies}
We briefly describe the strategies used in this paper.
\vskip0.1in
The starting point for the proof of Theorem \ref{thm1.1} is the argument of Hang and Wang in \cite{Hangwang}: we argue by contradiction and apply the concentration--compactness principle to obtain a sequence of functions $u_j$ with suitable convergence properties. The key step is to decompose $u_j$ into a sum of several Talenti bubbles and a remainder that is orthogonal to the tangent spaces of the Talenti bubbles and satisfies a spectral inequality (see Proposition \ref{prop2.3}). Through careful expansions, we derive delicate error estimates for this decomposition (see Lemmas \ref{lem2.4} and \ref{lem2.6}). Finally, we expand the Sobolev inequality using this decomposition and apply the error estimates to obtain a contradiction.
\vskip0.1in
To prove Theorem \ref{thm1.2}, we use Druet's localization theorem \cite[Theorem 4.9]{Heb1}, which states that the failure of the \textbf{sharp} Sobolev inequality \eqref{Sob3} is equivalent to its local failure near some point. Placing this local failure phenomenon around the support of an extremal measure yields the desired test function.
\vskip0.1in
The proof of Theorem \ref{thm1.4} relies on three ingredients. The first is a penalized higher-order Moser--Trudinger--Onofri inequality (see Lemma \ref{lem4.2}), which strengthens the classical inequality \eqref{MTO0} by an $L^2$-term. The second is a mass separation lemma (see Lemma \ref{lem4.3}), which provides a uniform description of the mass distribution of the measure $\frac{e^{nu}}
    {\int_{\mathbb S^n}e^{nu}\de \mu}\,\de \mu$. The third (see Lemma \ref{lem4.5}) provides estimates for an orthogonal decomposition of $H^s(\S^n)$ adapted to the mass distributions supplied by the second ingredient. To deduce Theorem \ref{thm1.4}, we first apply the mass separation lemma to reduce the left-hand side of \eqref{newmain2} to $\log \int_{U_i} e^{n u}\,\mathrm{d}\mu$ for some subset $U_i\subset \S^n$. The estimates for the orthogonal decomposition then reduce this expression further to $\log \int_{U_i} e^{n u_i}\,\mathrm{d}\mu$, where $u_i$ is an orthogonal component of $u$. We next apply the penalized higher-order Moser--Trudinger--Onofri inequality to each $u_i$. The additional $L^2$-term allows us to control the errors and ultimately deduce \eqref{newmain2}.
\vskip0.1in
    For the proof of Theorem \ref{thm1.5}, we construct test functions for general $m$ by suitably perturbing constant functions. This construction is effective for large $n$. For $m=1$, we also construct explicit test functions to cover the small-$n$ case.
\vskip0.1in
    To prove Theorems \ref{thm1.6} and \ref{thm5.1}, we construct a suitable $l^p$-partition of unity $\eta_j$, apply \eqref{Sob3} to each $\eta_ju$, and sum the resulting estimates. The key is to establish an effective lower bound for the term in \eqref{eq5.Phi-r} (see Proposition \ref{prop5.7}). This lower bound follows from a careful choice of partition functions (see Lemma \ref{lem5.4}) and a quantitative analysis of almost-orthogonality (see Lemmas \ref{lem5.2} and \ref{lem5.6}).
\vskip0.1in
    For Theorems \ref{thm1.7} and \ref{thm6.2}, we combine the mass separation lemma (Lemma \ref{lem4.3}) with a quantitative ``distribution of mass'' principle for the higher-order Moser--Trudinger--Onofri inequality \eqref{MTO0} (see Lemma \ref{lem6.1}). This type of principle originates in the work of
    Chen and Li \cite{Chenli} and has been widely used in the study of improved inequalities. However, standard approaches to this principle, such as truncation arguments and eigenspace decompositions, do not apply in our setting. We therefore develop a different method: we first reduce the estimate to a Riesz-type inequality using Green's representation and then apply a truncation argument to the resulting inequality.
\vskip0.1in
    Finally, to prove Corollaries \ref{cor1.8} and \ref{cor1.9}, respectively, we combine Theorems \ref{thm1.6} and \ref{thm1.7} with the elementary inequality $(x+y)^q\geq x^q+qx^{q-1}y$ for $x,y\geq 0$ and $q\geq 1$.    
\subsection{Structure of the paper}
In Section \ref{sec2}, we prove Theorem \ref{thm2.1}, which is equivalent to Theorem \ref{thm1.1}. Theorem \ref{thm1.2} is proved in Section \ref{sec3}. Section \ref{sec4} is devoted to the proofs of Theorems \ref{thm4.1} and \ref{thm4.6}, which imply Theorems \ref{thm1.4} and \ref{thm1.5}, respectively. In Section \ref{sec5} and Subsection \ref{sec6.1}, we establish Theorems \ref{thm5.1} and \ref{thm6.2}, which imply Theorems \ref{thm1.6} and \ref{thm1.7}, respectively. The proofs of Corollaries \ref{cor1.8} and \ref{cor1.9} are also given there. In Subsection \ref{sec6.2}, we generalize \cite[Theorem 1.1]{Chenzhang} to higher-order inequalities on $\S^n$ for all $n\geq 2$. Finally, in Section \ref{sec7}, we study improved inequalities under symmetry constraints and obtain various new results on sharpness, non-sharpness, and quantitative remainder estimates.

\medskip
\noindent\textbf{Note added after submission.}
The present manuscript was essentially completed in early
September 2026 and submitted to a journal on September 17.
The submitted version was later posted on the arXiv on
October 19. In the meantime, three closely related
preprints appeared: Wei, Yan and Ye \cite{WYY-Dolbeault}
on October 9, Frank, Jin, Wang and Xiong \cite{FJWX}
on October 14, and Wei, Yan and Ye \cite{WYY-Sobolev}
on October 16. Our work was carried out independently
of these preprints.

The case $m=1$ of Theorem \ref{thm1.1} overlaps with results
in \cite{FJWX,WYY-Sobolev}. In particular,
\cite{WYY-Sobolev} obtains the optimal $L^2$ coefficient
for $n\geq 4$, while \cite{FJWX} determines the exact
parameter range for the sharp inequality under
first-moment constraints, including a wider non-sharpness
range than that in Corollary \ref{cor1.3}.
The paper \cite{WYY-Dolbeault} establishes sharp
inequalities for antipodally even functions and certain
axially symmetric functions, including in dimension three.
Our Theorem \ref{thm1.1} also covers $m=2,3$, and our other
results address different questions. The statements
concerning novelty in the original manuscript reflect
our knowledge at the time of submission.
\medskip

\section{Sharp improved Sobolev inequalities under low-order moment constraints}\label{sec2}
Throughout this section, we assume $n\geq 3$ and let $\de w=\frac{\de\mu}{\omega_n}$ be the
normalized measure on $\S^n$. All norms are taken with respect to
$\de w$. We set
\[
 q=\frac{2n}{n-2},\qquad \alpha=\frac{2}{q}=\frac{n-2}{n},
 \qquad \kappa=\frac{n(n-2)}4,
\]
and equip $H^1(\S^n)$ with the norm
\[
 \norm{v}_H^2=E(v):=\int_{\S^n}
 \bigl(\abs{\nabla v}^2+\kappa v^2\bigr)\,\de w.
\]
The sharp Sobolev inequality on $\S^n$ reads (see, for example, \cite{Bec})
\begin{align}\label{sharpsob}
    E(v)\geq \kappa\norm{v}_{L^q(\S^n)}^2.
\end{align}
For $m\in\{1,2,3\}$, let
\[
 N_m=\begin{cases}
 2,&m=1,\\
 n+2,&m=2,\\
 2n+2,&m=3,
 \end{cases}
 \qquad c_m=N_m^{-1},\qquad
 \Theta_m=\Theta(m,\alpha,n)=N_m^{1-\alpha}.
\]
For
$z,\xi\in\S^n$, put $s=z\cdot\xi$ and define three polynomials
\begin{equation}
 H_{m,z}(\xi)=
 \begin{cases}
 \dfrac{1+s}{2},&m=1,\\[6pt]
 \dfrac{[1+(n+1)s]^2}{(n+2)^2},&m=2,\\[6pt]
 \dfrac{s^2(1+s)}{2},&m=3.
 \end{cases}
\end{equation}
They satisfy
\begin{equation}\label{eq2.certificate-properties}
 H_{m,z}\geq0,\qquad H_{m,z}(z)=1,
 \qquad \int_{\S^n}H_{m,z}\,\de w=c_m,
 \qquad H_{m,z}-c_m\in\PP_{m,n}.
\end{equation}

Our main result is the following, which is equivalent to Theorem \ref{thm1.1}.
\begin{theorem}\label{thm2.1}
For every $m\in\{1,2,3\}$, there is a constant $B=B(m,n)<\infty$ such
that every $u\in H^1(\S^n)$ satisfying \eqref{or1} obeys
\begin{equation}
 \norm{u}_{L^q(\S^n)}^2
 \leq \frac{1}{\kappa\Theta_m}
 \norm{\nabla u}_{L^2(\S^n)}^2
 +B\norm{u}_{L^2(\S^n)}^2.
\end{equation}
\end{theorem}
The proof of Theorem \ref{thm2.1} consists of a series of lemmas. The first lemma follows from \cite[Theorem 4.1]{Put}, which provides useful estimates for cubature formulas. We include a simple proof for completeness.
\begin{lemma}\label{lem2.2}
Suppose that
$\nu=\sum\limits_i w_i\delta_{z_i}\in\M^c_{m,n}$, where $w_i>0$ and
$\sum\limits_iw_i=1$. Then
\begin{equation}
 w_i\leq c_m\quad\text{for every }i,
 \qquad \sum_iw_i^\alpha\geq\Theta_m.
\end{equation}
Moreover, equality in the second inequality holds if and only if $\nu$ has exactly
$N_m$ atoms, all of weight $c_m=N_m^{-1}$, and the configurations of the points $z_i$ are, respectively, an antipodal pair when $m=1$, the vertices of a regular simplex when $m=2$, and the
vertices of an orthogonal cross-polytope when $m=3$.
\end{lemma}

\begin{proof}
From \eqref{eq2.certificate-properties} we have
\[
 w_i\leq\sum_k w_kH_{m,z_i}(z_k)
 =\int_{\S^n}H_{m,z_i}\,\de w=c_m.
\]
Since $\alpha-1<0$,
\[
 \sum_iw_i^\alpha=\sum_iw_iw_i^{\alpha-1}
 \geq c_m^{\alpha-1}\sum_iw_i
 =N_m^{1-\alpha}=\Theta_m.
\]
Equality forces every positive weight to equal $c_m$, and the converse is
immediate. Notice also that equality forces
$H_{m,z_i}(z_k)=0$ for $i\neq k$. Thus the equality configurations are,
respectively, an antipodal pair, the vertices of a regular simplex, and the
vertices of an orthogonal cross-polytope.
\end{proof}

Our next proposition provides a crucial refinement of the concentration--compactness argument used by Hang and Wang in \cite{Hangwang}. Before stating it, we recall some basic facts about Sobolev extremizers on $\S^n$. For $0<\varepsilon\leq1$ and
$z\in\S^n$, let
\begin{equation}\label{eq2.bubble}
 U_{\varepsilon,z}(\xi)
 =\left(\frac{2\varepsilon}
 {1+\varepsilon^2-(1-\varepsilon^2)\xi\cdot z}\right)^{\frac{n-2}{2}}.
\end{equation}
These are the normalized conformal Talenti bubbles. If we set
$A=-\Delta_{\S^n}+\kappa$, then
\begin{align}\label{eq2.bubble-identities}
 \norm{U_{\varepsilon,z}}_{L^q(\S^n)}=1,\qquad
 AU_{\varepsilon,z}=\kappa U_{\varepsilon,z}^{q-1},
 \qquad E(U_{\varepsilon,z})=\kappa.
\end{align}
Let $T_{\varepsilon,z}$ be a conformal isometry of both $H^1(\S^n)$ and $L^q(\S^n)$
such that $T_{\varepsilon,z}1=U_{\varepsilon,z}$. If $\mathcal H_k$
denotes the space of spherical harmonics of degree $k$, the tangent space to
the bubble manifold at $U=U_{\varepsilon,z}$ is
\begin{align}
 \mathcal T(U)=T_{\varepsilon,z}(\mathcal H_0\oplus\mathcal H_1).
\end{align}

\begin{proposition}\label{prop2.3}
Let $m\in\{1,2,3\}$, and let $u_j\in H^1(\S^n)$ with $u_j\geq0$ satisfy \eqref{or1},
\begin{equation}\label{eq2.critical-sequence}
 \norm{u_j}_{L^q(\S^n)}=1,\qquad \norm{u_j}_{L^2(\S^n)}\longrightarrow0,
 \qquad \limsup_{j\to\infty}E(u_j)\leq\kappa\Theta_m.
\end{equation}
After passing to a subsequence, there are distinct points
$z_1^*,\ldots,z_{N_m}^*\in\S^n$ such that
\begin{equation}\label{eq2.critical-measure}
 \abs{u_j}^q\de w\weakstarto
 \frac1{N_m}\sum_{i=1}^{N_m}\delta_{z_i^*},
 \qquad E(u_j)\longrightarrow\kappa\Theta_m.
\end{equation}
Moreover, one can write
\begin{equation}\label{eq2.modulation}
 u_j=\sum_{i=1}^{N_m}a_{i,j}U_{i,j}+r_j,
 \qquad U_{i,j}=U_{\varepsilon_{i,j},z_{i,j}},
\end{equation}
where
\begin{equation}\label{eq2.modulation-limits}
 a_{i,j}\longrightarrow N_m^{-\frac{1}{q}},\qquad
 \varepsilon_{i,j}\longrightarrow0,\qquad
 z_{i,j}\longrightarrow z_i^*,\qquad
 \norm{r_j}_H\longrightarrow0,
\end{equation}
and
\begin{equation}\label{eq2.orthogonality}
 r_j\perp_H\mathcal T(U_{i,j})\qquad(1\leq i\leq N_m).
\end{equation}
Furthermore, for every
$0<\gamma<\frac{4}{n+4}$ and all sufficiently large $j$,
\begin{equation}\label{eq2.spectral-gap}
 \norm{r_j}_H^2-(q-1)\kappa\sum_{i=1}^{N_m}
 \int_{\S^n}U_{i,j}^{q-2}r_j^2\,\de w
 \geq\gamma\norm{r_j}_H^2.
\end{equation}
\end{proposition}

\begin{proof}
Since $u_j\weakto0$ in $H^1(\S^n)$, the concentration--compactness principle (see for example
\cite[Theorem~2.1]{Hangwang}) yields that, up to a subsequence,
\[
 \abs{u_j}^q\de w\weakstarto\nu=\sum_\ell\nu_\ell\delta_{x_\ell},
 \qquad \abs{\nabla u_j}^2\de w\weakstarto\mu,
 \qquad
 \mu\geq\kappa\sum_\ell\nu_\ell^\alpha\delta_{x_\ell}.
\]
The moment constraints pass to the limit, so $\nu\in\M^c_{m,n}$.
Lemma~\ref{lem2.2} and \eqref{eq2.critical-sequence} therefore give
\[
 \kappa\Theta_m\geq\mu(\S^n)
 \geq\kappa\sum_\ell\nu_\ell^\alpha
 \geq\kappa\Theta_m,
\]
Thus every inequality is an equality. Lemma~\ref{lem2.2} then yields
\eqref{eq2.critical-measure}.

Choose fixed disjoint neighborhoods of the points $z_i^*$ and cut off
$u_j$ in these neighborhoods. Write $u_j=\sum\limits_{i=1}^{N_m}\chi_iu_j+\Tilde{u}_j$, where $\chi_i$ is a cutoff function near $z_i^*$. The convergence properties imply that $\norm{\Tilde{u}_j}_H\rightarrow 0$, and each $\chi_iu_j$ has $L^q$-mass tending to $N_m^{-1}$ and
energy tending to $\kappa N_m^{-\alpha}$. Thus the Sobolev deficit of $\chi_iu_j$ tends to zero. By the spherical form of the Bianchi--Egnell stability
theorem \cite{Bia}, we can replace $\chi_iu_j$ with $a_{i,j}^0U_{i,j}^0$, where $a_{i,j}^0\rightarrow N_m^{-\frac{1}{q}}$, $U_{i,j}^0=U_{\varepsilon^0_{i,j},z^0_{i,j}}$ is a Talenti bubble, $\varepsilon^0_{i,j}\rightarrow0$, and $z_{i,j}^0\rightarrow z_i^*$. Moreover, we have $r_j^0:=u_j-\sum\limits_{i=1}^{N_m}a_{i,j}^0U_{i,j}^0\rightarrow 0$ in $H^1(\S^n)$.

We next modify this decomposition to obtain the orthogonality condition
\eqref{eq2.orthogonality}.
Write $T_{i,j}^0=T_{\varepsilon_{i,j}^0,z_{i,j}^0}$. We use normalized coordinates on the bubble manifold. For
$\vartheta=(s,b)\in\R\times\R^{n+1}$ near the origin, set
\begin{equation}
 W(s,b)(\xi)=(1+s)
 \left(\frac{1-\abs{b}^2}
 {1+\abs{b}^2-2b\cdot\xi}\right)^{\frac{n-2}{2}}.
\end{equation}
Thus $W(0,0)=1$, and the derivatives of $W$ at the origin span
$\mathcal H_0\oplus\mathcal H_1$. Consequently, for
$\boldsymbol\vartheta=(\vartheta_1,\ldots,\vartheta_{N_m})\in
\R^{d}$, $d=N_m(n+2)$, the map
\begin{equation}
 \Psi_j(\boldsymbol\vartheta)
 =\sum_{i=1}^{N_m}a_{i,j}^0T_{i,j}^0W(\vartheta_i)
\end{equation}
parametrizes a neighborhood of the initial
multi-bubble $\sum\limits_{i=1}^{N_m}a_{i,j}^0U_{i,j}^0$.

Index the coordinates of $\boldsymbol\vartheta$ by
$\beta=(i,k)$, $1\leq i\leq N_m$, $0\leq k\leq n+1$, and put
\[
 \phi_{j,\beta}(\boldsymbol\vartheta)
 =\partial_{\vartheta_\beta}\Psi_j(\boldsymbol\vartheta).
\]
For a fixed $i$, the vectors
$\{\phi_{j,(i,k)}(\boldsymbol\vartheta):0\leq k\leq n+1\}$
span the tangent space of the $i$th bubble factor. Define the
finite-dimensional map
\begin{equation}
 F_j(\boldsymbol\vartheta)
 =\left(
 \left\langle u_j-\Psi_j(\boldsymbol\vartheta),
 \phi_{j,\beta}(\boldsymbol\vartheta)\right\rangle_H
 \right)_{\beta}\in\R^d.
\end{equation}
Hence $F_j(\boldsymbol\vartheta)=0$ is equivalent to the orthogonality
conditions in \eqref{eq2.orthogonality}.

By conformal invariance and our choice of coordinates, the $C^3$-norm of $\Psi_j$ is uniformly bounded in a fixed neighborhood of $0$. We therefore compute the derivative at the initial parameter
$\boldsymbol\vartheta=0$:
\begin{equation}
 \partial_{\vartheta_\alpha}F_{j,\beta}(0)
 =-\left\langle\phi_{j,\alpha}(0),
 \phi_{j,\beta}(0)\right\rangle_H
 +\left\langle r^0_j,
 \partial_{\vartheta_\alpha}\phi_{j,\beta}(0)
 \right\rangle_H=-(G_j)_{\alpha\beta}+o(1),
\end{equation}
where $(G_j)_{\alpha\beta}
 =\left\langle\phi_{j,\alpha}(0),
 \phi_{j,\beta}(0)\right\rangle_H$.
The diagonal blocks of $G_j$ converge to positive multiples of the fixed
single-bubble Gram matrix
\[
 G_*=\left(
 \left\langle\partial_{\vartheta_k}W(0),
 \partial_{\vartheta_\ell}W(0)\right\rangle_H
 \right)_{0\leq k,\ell\leq n+1},
\]
which is positive definite. Since $\varepsilon^0_{i,j}\rightarrow0$, $z_{i,j}^0\rightarrow z_i^*$, and $z_i^*\neq z_j^*$ for $i\neq j$, all
off-diagonal blocks tend to zero. Therefore,
\begin{equation}
 G_j\longrightarrow
 N_m^{-\frac{2}{q}}G_*\oplus\cdots\oplus N_m^{-\frac{2}{q}}G_*,
\end{equation}
so $DF_j(0)$ is invertible for large $j$ and
$\norm{DF_j(0)^{-1}}\leq C$ uniformly in $j$.

Since the $C^3$-norm of $\Psi_j$ is uniformly bounded in a fixed neighborhood of $0$ and $|F_j(0)|\rightarrow 0$, we may use the local inverse function theorem in the following form: if
$F:\R^d\to\R^d$ is $C^1$, $DF(0)^{-1}$ is bounded by $K$, $DF$ is sufficiently close to $DF(0)$ in a fixed ball, and $|F(0)|$ is small, then there is a unique zero $\vartheta$ of $F$ with $\abs{\vartheta}\leq2K\abs{F(0)}$. Applying this theorem to $F_j$ for large $j$, we find a unique
$\boldsymbol\vartheta_j$ such that
$F_j(\boldsymbol\vartheta_j)=0$ and
$\abs{\boldsymbol\vartheta_j}\rightarrow 0$. Set $a_{i,j}U_{i,j}=a_{i,j}^0T_{i,j}^0W(\vartheta_{i,j})$ and $r_j=u_j-\sum\limits_{i=1}^{N_m}a_{i,j}U_{i,j}$. Then \eqref{eq2.modulation-limits} and \eqref{eq2.orthogonality} follow.

It remains to prove \eqref{eq2.spectral-gap}. If not, we could
normalize $\norm{r_j}_H=1$ and, after pulling back by each bubble, assume
\[
 T_{\varepsilon_{i,j},z_{i,j}}^{-1}r_j\weakto R_i
 \quad\text{in }H^1(\S^n),\qquad
 T_{\varepsilon_{i,j},z_{i,j}}^{-1}r_j\longrightarrow R_i
 \quad\text{in }L^2(\S^n).
\]
By \eqref{eq2.orthogonality},
$R_i\perp_H(\mathcal H_0\oplus\mathcal H_1)$. Expanding $\norm*{r_j-\sum\limits^{N_m}_{i=1}T_{\varepsilon_{i,j},z_{i,j}}R_i}_H^2\geq 0$, we obtain
$\sum\limits^{N_m}_{i=1}\norm{R_i}_H^2\leq1$. Since the first spherical eigenvalue outside
$\mathcal H_0\oplus\mathcal H_1$ is $2(n+1)$,
\[
 (q-1)\kappa\sum_{i=1}^{N_m}\norm{R_i}_{L^2(\S^n)}^2
 \leq\frac{(q-1)\kappa}{\kappa+2(n+1)}
 \sum_{i=1}^{N_m}\norm{R_i}_H^2\leq\frac{n}{n+4}.
\]
Finally, conformal invariance gives
\[
 \int_{\S^n}U_{i,j}^{q-2}r_j^2\,\de w
 =\norm{T_{\varepsilon_{i,j},z_{i,j}}^{-1}r_j}_{L^2(\S^n)}^2\rightarrow \norm{R_i}_{L^2(\S^n)}^2.
\]
Thus the lower limit of the left-hand side of
\eqref{eq2.spectral-gap} is at least $\frac{4}{(n+4)}$, which yields a contradiction. This completes the proof.
\end{proof}
In what follows, we work along a sequence with the decomposition
\eqref{eq2.modulation}--\eqref{eq2.orthogonality}. We henceforth assume that $j$ is sufficiently large and suppress the index
$j$. Set
\[
 V=\sum_{i=1}^{N_m}a_iU_i,\qquad t=\norm{r}_H,\qquad
 M=\sum_{i=1}^{N_m}\ell_n(\varepsilon_i),
\]
where
\begin{equation}
 \ell_n(\varepsilon)=
 \begin{cases}
 \varepsilon,&n=3,\\
 \varepsilon^2\log(1/\varepsilon),&n=4,\\
 \varepsilon^2,&n\geq5.
 \end{cases}
\end{equation}

Our next two lemmas establish further estimates for this decomposition.
\begin{lemma}\label{lem2.4}
Assume in addition that $\norm{u}_{L^q(\S^n)}=1$ and that $u=V+r$ satisfies \eqref{or1}. Then there
exist positive constants $b$ and $C$ such that
\begin{equation}\label{eq2.amplitude-estimate}
 \left|\sum_{i=1}^{N_m}a_i^q-1\right|
 +\max_{1\leq i\leq N_m}\left|a_i^q-c_m\right|
 \leq C(M+t^2),
\end{equation}
and
\begin{equation}\label{eq2.local-deficit}
 E(u)-\kappa\Theta_m\norm{u}_{L^q(\S^n)}^2
 \geq b\norm{r}_H^2-CM.
\end{equation}
\end{lemma}

\begin{proof}
With $q'=q/(q-1)$, set
\[
 I_{ik}=\int_{\S^n}U_i^{q-1}U_k\,\de w,\qquad
 \delta=\left\|V^{q-1}-\sum^{N_m}_{i=1} a_i^{q-1}U_i^{q-1}
 \right\|_{L^{q'}(\S^n)}.
\]
The explicit decay in \eqref{eq2.bubble} gives
\begin{equation}\label{eq2.interactions}
 \sum_{i<k}I_{ik}\leq CM,\qquad \delta^2\leq CM,\qquad
 \left\|V^{q-2}-\sum^{N_m}_{i=1} a_i^{q-2}U_i^{q-2}
 \right\|_{L^{\frac{n}{2}}(\S^n)}=o(1).
\end{equation}
Moreover, uniformly for $f$ in a bounded
subset of $C^2(\S^n)$, we have
\begin{equation}\label{eq2.moving-test}
 \int_{\S^n}fU_i^q\,\de w=f(z_i)+O_f(\varepsilon_i^2),
 \qquad
 \norm{(f-f(z_i))U_i^{q-1}}_{L^{q'}(\S^n)}\leq C_f\varepsilon_i.
\end{equation}
The orthogonality condition \eqref{eq2.orthogonality} yields
\begin{equation}\label{eq2.amplitude-orthogonality}
 \int_{\S^n}U_i^{q-1}r\,\de w=0.
\end{equation}
For $f$ in a bounded
subset of $C^2(\S^n)$, the first-order expansion of
$\abs{V+r}^q$ gives
\[
 \int_{\S^n}f\abs{V+r}^q\,\de w
 =\int_{\S^n}fV^q\,\de w
 +q\int_{\S^n}fV^{q-1}r\,\de w+O_f(t^2).
\]
Combining \eqref{eq2.interactions}, \eqref{eq2.moving-test}, and
\eqref{eq2.amplitude-orthogonality} with a standard expansion argument, we obtain
\[
 \left|\int_{\S^n}fV^{q-1}r\,\de w\right|
 \leq C_f\left(\delta+\sum^{N_m}_{i=1}\varepsilon_i\right)t,
 \qquad
 \int_{\S^n}fV^q\,\de w
 =\sum^{N_m}_{i=1} a_i^qf(z_i)+O_f(M).
\]
Since $\delta^2+\sum\limits^{N_m}_{i=1}\varepsilon_i^2\leq CM$, Young's inequality yields
\begin{equation}
 \int_{\S^n}f\abs{u}^q\,\de w
 =\sum^{N_m}_{i=1} a_i^qf(z_i)+O_f(M+t^2).
\end{equation}
Taking $f=1$ gives
\begin{align}\label{qd}
    T:=\sum^{N_m}_{i=1} a_i^q=1+O(M+t^2).
\end{align}

Next, fix $i$ and take
$f=H_{m,z_i}$. Equations \eqref{eq2.certificate-properties} and \eqref{or1} give
\[
 c_m=\int_{\S^n}H_{m,z_i}\abs{u}^q\,\de w
 =\sum^{N_m}_{k=1} a_k^qH_{m,z_i}(z_k)+O(M+t^2).
\]
The nonnegativity of $H_{m,z_i}$ and the identity $H_{m,z_i}(z_i)=1$ imply
$a_i^q-c_m\leq C(M+t^2)$. From \eqref{qd} and $N_mc_m=1$, we also obtain $a_i^q-c_m\geq -C(M+t^2)$. This proves \eqref{eq2.amplitude-estimate}.

We turn to the estimate \eqref{eq2.local-deficit}. From \eqref{eq2.bubble-identities} we have
\begin{equation}\label{eq2.energy-expansion}
 E(u)=\kappa\sum^{N_m}_{i=1} a_i^2
 +2\kappa\sum_{i<k}a_ia_kI_{ik}+t^2.
\end{equation}
The second-order expansion yields
\begin{equation}
 \int_{\S^n}\abs{V+r}^q\,\de w
 =\int_{\S^n}V^q\,\de w
 +q\int_{\S^n}V^{q-1}r\,\de w
 +\frac{q(q-1)}2\int_{\S^n}V^{q-2}r^2\,\de w+o(t^2).
\end{equation}
Using the concavity of $s\mapsto s^\alpha$ together with
\eqref{eq2.interactions} and \eqref{eq2.energy-expansion}, we obtain
\begin{align}
 E(u)-\kappa\Theta_m\norm{u}_{L^q(\S^n)}^2
 \geq{}&\kappa\left(\sum_{i=1}^{N_m} a_i^2-\Theta_mT^\alpha\right)+t^2\nonumber\\
 &-(q-1)\kappa\Theta_mT^{\alpha-1}
 \int_{\S^n}V^{q-2}r^2\,\de w
 -C(M+\delta t)+o(t^2).
\end{align}
To estimate the first term on the right, set $x_i=a_i^q$ and $\bar x=\frac{T}{N_m}$. Combining
\eqref{eq2.amplitude-estimate} and Taylor's
formula with $\sum\limits^{N_m}_i(x_i-\bar x)=0$ gives
\begin{equation}\label{eq2.discrete-holder}
\sum_{i=1}^{N_m} a_i^2-\Theta_mT^\alpha=\sum_{i=1}^{N_m}x_i^\alpha-N_m\bar x^\alpha\geq-C\sum_{i=1}^{N_m}(x_i-\bar x)^2\geq -CM+o(t^2).
\end{equation}
On the other hand, \eqref{eq2.amplitude-estimate} gives
\begin{equation}\label{eq2.coefficient-matching}
 \Theta_mT^{\alpha-1}a_i^{q-2}=1+O(M+t^2).
\end{equation}
It follows from \eqref{eq2.discrete-holder},
\eqref{eq2.interactions} and \eqref{eq2.coefficient-matching} that
\begin{align}
    E(u)-\kappa\Theta_m\norm{u}_{L^q(\S^n)}^2\geq t^2-(q-1)\kappa\sum_{i=1}^{N_m}
 \int_{\S^n}U_i^{q-2}r^2\,\de w-C(M+\delta t)+o(t^2).\nonumber
\end{align}
By applying \eqref{eq2.spectral-gap}, \eqref{eq2.interactions}, and Young's inequality
$\delta t\leq\eta t^2+C_\eta \delta^2$ for $\eta$ small, we obtain \eqref{eq2.local-deficit}.
\end{proof}
\begin{remark}\label{rmk2.5}
    We now explain why the assumption $m\in\{1,2,3\}$ is crucial in our proof. By Lemma \ref{lem2.2}, every extremal measure $\nu$ has exactly $N_m$ atoms, all of weight $c_m=N_m^{-1}$. This property makes the first term in the expansion of $E(u)-\kappa\Theta_m\norm{u}_{L^q(\S^n)}^2$ small (see \eqref{eq2.discrete-holder}) and
    makes the quadratic term exactly the left-hand side of \eqref{eq2.spectral-gap}, yielding $ct^2$ as a lower bound. When $m\geq 4$, the configurations of extremal measures are generally difficult to determine (see, for example, \cite{de,Del,Put}); if the weights of $\nu$ are not equal, one may not be able to apply \eqref{eq2.spectral-gap}.
\end{remark}

\begin{lemma}\label{lem2.6}
For every fixed sufficiently small constant $\rho>0$, we have
\begin{equation}\label{eq2.l2-absorption}
 M\leq C_\rho\norm{u}_{L^2(\S^n)}^2
 +C_0\eta_n(\rho)\norm{r}_H^2,
 \qquad
 \eta_n(\rho)=
 \begin{cases}
 \rho,&n=3,\\
 \rho^2,&n\geq4,
 \end{cases}
\end{equation}
where $C_0$ is independent of $\rho$.
\end{lemma}
\begin{proof}
Since $j$ is sufficiently large, the balls $B_{i}=B(z_{i},\rho)$ are disjoint, and $a_i\geq \frac{N_m^{-\frac{1}{q}}}{2}$. We then have
\begin{equation}\label{eq2.local-mass-upper}
 \sum_{i=1}^{N_m}\int_{B_{i}}a_i^2U_{i}^2\,\de w
 \leq \sum_{i=1}^{N_m}\int_{B_i}V^2\de w\leq 2\norm{u}_{L^2(\S^n)}^2
 +2\sum_{i=1}^{N_m}\int_{B_{i}}r^2\,\de w.
\end{equation}
H\"older's inequality and the sharp Sobolev inequality \eqref{sharpsob} yield
\[
 \sum_{i=1}^{N_m}\int_{B_{i}}r^2\,\de w
 \leq C\rho^2\norm{r}_H^2.
\]
Direct radial integration in \eqref{eq2.bubble} gives
\[
 \int_{B_{i}}U_{i}^2\,\de w\geq
 \begin{cases}
 c\rho\varepsilon_{i},&n=3,\\
 -c\varepsilon_{i}^2\log(\varepsilon_{i}),&n=4,\\
 c\varepsilon_{i}^2,&n\geq5.
 \end{cases}
\]
Inserting these two estimates into \eqref{eq2.local-mass-upper} yields
\eqref{eq2.l2-absorption}.
\end{proof}
Now we are ready to give the proof of Theorem \ref{thm2.1}.
\begin{proof}[Proof of Theorem~\ref{thm2.1}]
Fix $m\in\{1,2,3\}$ and define
\[
 D_m(u)=E(u)-\kappa\Theta_m\norm{u}_{L^q(\S^n)}^2.
\]
It suffices to prove that $D_m(u)\geq-B_{m,n}\norm{u}_{L^2}^2$ for some
finite constant $B_{m,n}$. Otherwise, we can find nonnegative functions
$u_j\in H^1(\S^n)$ satisfying \eqref{or1} such that
\begin{equation}\label{eq2.contradiction-sequence}
 \norm{u_j}_{L^q(\S^n)}=1,\qquad
 D_m(u_j)<-j\norm{u_j}_{L^2(\S^n)}^2.
\end{equation}
Consequently,
\[
 E(u_j)+j\norm{u_j}_{L^2(\S^n)}^2<\kappa\Theta_m,
\]
so \eqref{eq2.critical-sequence} holds. Proposition~\ref{prop2.3}
and Lemma~\ref{lem2.4} give
\[
 D_m(u_j)\geq b\norm{r_j}_H^2-CM_j.
\]
Choose $\rho>0$ so small that
$CC_0\eta_n(\rho)<b/2$. Lemma~\ref{lem2.6} then gives, for all
sufficiently large $j$,
\[
 D_m(u_j)\geq\frac b2\norm{r_j}_H^2
 -CC_\rho\norm{u_j}_{L^2(\S^n)}^2
 \geq-CC_\rho\norm{u_j}_{L^2(\S^n)}^2,
\]
contradicting \eqref{eq2.contradiction-sequence}. This completes the proof.
\end{proof}

\section{A failure principle for sharp improved Sobolev inequalities}\label{sec3}
In this section, we prove Theorem \ref{thm1.2}. We first recall Druet's localization theorem \cite[Theorem 4.9]{Heb1}, whose main idea originates in Druet's work \cite{Dru}.

\begin{theorem}\label{thm3.1}
    Let $(M,g)$ be a smooth closed manifold of dimension $n\geq 2$. Suppose that, for some $1\leq p<n$, the \textbf{sharp} Sobolev inequality \eqref{Sob3} is locally valid in the sense that each $x\in M$ possesses an open neighborhood $\Omega_x$ such that, for every $u\in C_c^\infty(\Omega_x)$,
    \begin{align}
        \norm{u}_{L^{p^*}(M)}^p
     \leq S_{n,p}^p\norm{\nabla u}_{L^p(M)}^p
     + B_x\norm{u}_{L^{p}(M)}^p\nonumber
    \end{align}
    for some constant $B_x$ independent of $u$. Then the \textbf{sharp} Sobolev inequality \eqref{Sob3} holds globally on $M$ for some constant $B(M,p)$.
\end{theorem}
The following corollary is crucial to our construction.
\begin{corollary}\label{cor3.2}
Assume $1\leq p<n$ and that the \textbf{sharp} Sobolev inequality \eqref{Sob3} fails on $\mathbb S^n$. Then, for any $z\in \S^n$, there exist
a sequence $r_k\rightarrow0^+$ and nonnegative functions $w_k\in W^{1,p}(\mathbb S^n)$ such that
$\operatorname{supp}w_k\subset B(z,r_k)$, each $w_k$ is radially symmetric with respect to $z$, and
\[
\frac{
\|w_k\|_{L^{p^*}(\S^n)}^{p}
-
S_{n,p}^{p}\|\nabla w_k\|_{L^p(\S^n)}^{p}
}
{\|w_k\|_{L^p(\S^n)}^{p}}
\rightarrow+\infty .
\]
\end{corollary}

\begin{proof}
By Theorem \ref{thm3.1} and the rotation invariance of $\S^n$, for any $z\in \S^n$, we can find a sequence $r_k\rightarrow0^+$ and functions $v_k\in C_c^\infty(\mathbb S^n)$ such that
$\operatorname{supp}v_k\subset B(z,r_k)$ and
\[
\frac{
\|v_k\|_{L^{p^*}(\S^n)}^{p}
-
S_{n,p}^{p}\|\nabla v_k\|_{L^p(\S^n)}^{p}
}
{\|v_k\|_{L^p(\S^n)}^{p}}
\rightarrow+\infty .
\]
By the spherical P\'olya--Szeg\H{o} inequality, it suffices to take $w_k=|v_k|^*$, where $w_k$ is the spherical decreasing rearrangement of $|v_k|$ with respect to $z$.
\end{proof}
We now prove Theorem \ref{thm1.2}.
\begin{proof}[Proof of Theorem \ref{thm1.2}]
It suffices to find a sequence of test functions $u_k\in W^{1,p}(\mathbb S^n)$ satisfying \eqref{or1} such that
\begin{align}
    \frac{\|u_k\|_{L^{p^*}(\S^n)}^{p}-
\frac{S_{n,p}^{p}}{\Theta\left(m,\frac{n-p}{n},n\right)}
\|\nabla u_k\|_{L^p(\S^n)}^{p}}{
\|u_k\|_{L^p(\S^n)}^{p}}\rightarrow+\infty.\nonumber
\end{align}

By \cite[Corollary 3.2]{Put}, there exists a finitely supported measure $\nu\in \M^c_{m,n}$ such that $\norm{\nu}_{\frac{n-p}{n}}=\Theta(m,\frac{n-p}{n},n)$. Write $\nu=\sum\limits_{i=1}^{N}\nu_i\delta_{\xi_i}$.

Fix a point $z\in\S^n$ and let $w_{k}(y)=\phi_k(d(y,z))$ be the sequence of functions given by Corollary \ref{cor3.2}. Set $w_{i,k}(y)=\phi_k(d(y,\xi_i))$ for $1\leq i\leq N$. For sufficiently large $k$, the supports of $w_{i,k}$ are mutually disjoint. For such $k$, define
\[
u_k
=
\sum_{i=1}^{N}
\nu_i^{1/p^*}w_{i,k},
\]
We claim that these $u_k$ are the required test functions.

We first verify the moment constraints. Let $\mathcal{H}_k$ denote the space of spherical harmonics of degree
$k$ on $\mathbb S^{n}$. We have the orthogonal decomposition $\PP_{m,n}
=
\bigoplus_{k=1}^{m}\mathcal{H}_k$ (see \cite{Stein} for details). Since $w_k$ is radial, the
Funk--Hecke formula gives, for every
$Y\in\mathcal H_\ell$,
$1\leq\ell\leq m$,
\[
\int_{\mathbb S^n}Y(x)w_{i,k}^{p^*}(x)\,\mathrm{d}\mu(x)
=
b_{\ell,k}Y(\xi_i)
\]
for some constant $b_{\ell,k}$ independent of $i$. Therefore,
\begin{align}
\int_{\mathbb S^n}Y|u_k|^{p^*}\,\mathrm{d}\mu=
b_{\ell,k}\sum_{i=1}^{N}\nu_iY(\xi_i)=
b_{\ell,k}
\int_{\mathbb S^n}Y\,\de\nu
=0.\nonumber
\end{align}
Hence \eqref{or1} holds for $u_k$.

Since the supports of $w_{i,k}$ are disjoint, our choice of $\nu$ gives
\begin{align}
    \|u_k\|_{L^{p^*}(\S^n)}^{p^*}
=
\sum_{i=1}^{N}
\nu_i
\|w_k\|_{L^{p^*}(\S^n)}^{p^*}
=
\|w_k\|_{L^{p^*}(\S^n)}^{p^*}.\nonumber
\end{align}
Similarly, we have
\begin{align}
    \|\nabla u_k\|_{L^p(\S^n)}^p
=
\sum_{i=1}^{N}\nu_i^{p/p^*}
\|\nabla w_k\|_{L^p(\S^n)}^p=\Theta\left(m,\frac{n-p}{n},n\right)\|\nabla w_k\|_{L^p(\S^n)}^p,\nonumber
\end{align}
and
\begin{align}
    \|u_k\|_{L^p(\S^n)}^p
=
\sum_{i=1}^{N}\nu_i^{p/p^*}
\|w_k\|_{L^p(\S^n)}^p=\Theta\left(m,\frac{n-p}{n},n\right)\| w_k\|_{L^p(\S^n)}^p.\nonumber
\end{align}
Therefore,
\begin{align}
\frac{\|u_k\|_{L^{p^*}(\S^n)}^{p}-
\frac{S_{n,p}^{p}}{\Theta\left(m,\frac{n-p}{n},n\right)}
\|\nabla u_k\|_{L^p(\S^n)}^{p}}{
\|u_k\|_{L^p(\S^n)}^{p}}=
\frac1{\Theta\left(m,\frac{n-p}{n},n\right)}
\frac{
\|w_k\|_{p^*}^{p}
-
S_{n,p}^{p}\|\nabla w_k\|_{p}^{p}
}
{\|w_k\|_{p}^{p}}.\nonumber
\end{align}
Since the right-hand side tends to $+\infty$ by Corollary \ref{cor3.2}, this completes the proof.
\end{proof}

\section{Sharp improved Moser--Trudinger--Onofri inequalities under general moment constraints}\label{sec4}
In this section, we prove Theorems \ref{thm1.4} and \ref{thm1.5}. The first subsection deals with Theorem \ref{thm1.4}, and the second is devoted to Theorem \ref{thm1.5}.

Throughout this section, we write $\de w=\frac{\de\mu}{\omega_n}$ for the
normalized surface measure on the sphere \(\mathbb S^n\), and set $L=-\Delta_{\mathbb S^n}$ for notational simplicity. All norms are taken with respect to
$\de w$. Let \(\mathcal H_k\) be the space of spherical harmonics of degree \(k\). Then
\[
    L\big|_{\mathcal H_k}=\lambda_k\operatorname{Id},
    \qquad
    \lambda_k=k(k+n-1).
\]
For $u\in H^{\frac{n}{2}}(\S^n)$, define
\[
    E_n(u)
    :=
    \left\langle u,L^{n/2}u\right\rangle_{L^2(\de w)}
    =
    \sum_{k=1}^\infty
    \lambda_k^{\frac{n}{2}}\|u_k\|_{L^2(\S^n)}^2,
\]
where the operator $L^\frac{n}{2}$ is defined by
\[
    L^{\frac{n}{2}} u=\sum_{k=1}^\infty \lambda_k^\frac{n}{2} u_k.
\]
Here $u=\sum\limits_{k=0}^\infty u_k$ is the spherical harmonic decomposition. For $u,v\in H^{\frac{n}{2}}(\S^n)$, we also define
\begin{align}
    B_n(u,v):=\frac{1}{4}(E_n(u+v)-E_n(u-v))=\left\langle L^\frac{n}{4}u,L^{\frac{n}{4}}v\right\rangle_{L^2(\de w)}.\nonumber
\end{align}
If $\int_{\S^n} u\,\mathrm{d}w = 0$, that is, if $u_0=0$, then the Poincar\'e inequality gives $E_n(u)\geq n^\frac{n}{2}\norm{u}_{L^2(\S^n)}^2$.

Throughout what follows, we fix $m\geq 1$ and $n\geq 2$, and set
\begin{align}\label{yu}
    N_*=\Theta(m,0,n),\qquad A_{m,n}=\frac{n}{2(n-1)!\Theta(m,0,n)}.
\end{align}

\subsection{The sharpness result}

The following result is equivalent to Theorem \ref{thm1.4}.
\begin{theorem}\label{thm4.1}
There exists a positive constant $C=C(m,n)$ such that for any $u\in H^\frac{n}{2}(\S^n)$ satisfying
$\int_{\S^n} u\,\mathrm{d}w = 0$ and \eqref{or2},
\begin{align}
    \log \int_{\S^n} e^{n u}\,\mathrm{d}w
    \leq A_{m,n}E_n(u)
    + C.
\end{align}
\end{theorem}

The proof of Theorem \ref{thm4.1} has three ingredients. The first is the following penalized Adams inequality and its consequence, a penalized higher-order Moser--Trudinger--Onofri inequality. These strengthen the classical inequalities \eqref{MT3} and \eqref{MTO0}, respectively, in the case $s=\frac{n}{2}$.

\begin{lemma}\label{lem4.2}
Assume \(0\leq\kappa<n^{n/2}\) and define
\[
    Q_\kappa(v):=E_n(v)-\kappa\|v\|_{L^2(\S^n)}^2
\]
for every nonzero $v\in H^\frac{n}{2}(\S^n)$ satisfying
$\int_{\S^n} v\,\mathrm{d}w = 0$. Then there exists a positive constant $C(n,\kappa)$ such that
\begin{align}\label{penal1}
    \int_{\mathbb S^n}
    \exp\left(
       \frac{n!}{2}\frac{v^2}{Q_\kappa(v)}
    \right)\,\de w
    \leq C(n,\kappa).
\end{align}
Consequently,
\begin{equation}\label{penal2}
    \log\int_{\mathbb S^n}e^{nv}\,\de w
    \leq
    \frac{n}{2(n-1)!}Q_\kappa(v)+\log C(n,\kappa).
\end{equation}
\end{lemma}
Similar penalized Moser--Trudinger-type inequalities have been widely studied in the Euclidean setting; see, for example, \cite{Adi,Tin,Wang}. To the best of our knowledge, the inequalities \eqref{penal1} and \eqref{penal2} have not
appeared in this form. The key idea is to reduce \eqref{penal1} to a special case of the abstract Adams inequality of Fontana and Morpurgo \cite[Theorem 1]{Fon2}. This reduction uses some basic facts from the theory of pseudodifferential operators. We refer the reader to \cite[Sections 18.1--18.2]{Hor2} and \cite[Chapter 1]{Shu} for background.
\begin{proof}
Let $\Pi_0$ be the $L^2(\de w)$-orthogonal projection onto the constants. We define an operator
\begin{align}
    \mathcal B_\kappa
    =\Tilde{L}^\frac{n}{2}-\kappa(I-\Pi_0)=
    L^{n/2}-\kappa(I-\Pi_0)+\Pi_0,\nonumber
\end{align}
where $\Tilde{L}:=L+\Pi_0$. Note that $\mathcal B_\kappa$ acts as the identity on the constant
functions, while on $\mathcal H_k$, $k\geq1$, it acts by $\lambda_k^{n/2}-\kappa$. Since $\lambda_1=n$ and $\kappa<n^{n/2}$, all eigenvalues are strictly
positive. Hence $\mathcal B_\kappa$ is positive definite, self-adjoint, and invertible.

Since \(L\) is a second-order elliptic differential
operator and \(\Pi_0\) has finite rank,
\(\widetilde L=L+\Pi_0\) is a positive classical elliptic
pseudodifferential operator of order $2$, with principal symbol
\(|\zeta|_g^2\). Seeley's theorem on complex powers of elliptic operators (see, for example, \cite{See} and \cite[Chapter 2]{Shu}) implies that $\Tilde{L}^\frac{n}{2}$ is a classical elliptic
pseudodifferential operator of order \(n\), with principal symbol $|\zeta|_g^n$. It follows that \(\mathcal B_\kappa\) is also a classical elliptic
pseudodifferential operator of order \(n\), with principal symbol $|\zeta|_g^n$.

Applying Seeley's theorem on complex powers of elliptic operators once more (see, for example, \cite{See} and \cite[Chapter 2]{Shu}) to $\mathcal B_\kappa$, we find that $T_\kappa:=\mathcal B_\kappa^{-1/2}$ is a classical self-adjoint elliptic
pseudodifferential operator of order $-\frac{n}{2}$, with principal symbol $|\zeta|_g^{-\frac{n}{2}}$.

By \cite[Theorem 5.2.1]{Hor1}, we may take $K_\kappa(x,y)$ to be the Schwartz kernel of $T_\kappa$ relative to $\de w$. Since
$T_\kappa$ is a classical pseudodifferential operator, $K_\kappa(x,y)$ is smooth away from the diagonal. Moreover, since the principal symbol is $|\zeta|_g^{-\frac{n}{2}}$, a local kernel expansion near the diagonal shows that, for some $\delta>0$,
\begin{align}
    K_\kappa(x,y)
    =
    \omega_n(2\pi)^{-\frac{n}{2}} d_g(x,y)^{-\frac{n}{2}}
    +
    O(d_g(x,y)^{-\frac{n}{2}+\delta}),
    \qquad y\rightarrow x ,\nonumber
\end{align}
uniformly with respect to $x\in \S^n$ (the factor \(\omega_n\) appears because the kernel is taken with
respect to the normalized measure \(dw\)).

To apply the abstract Adams inequality \cite[Theorem 1]{Fon2} to $T_\kappa$, it suffices to check that conditions (13) and (14) therein hold for the kernel $K_\kappa$.

For fixed $x\in\mathbb S^n$, let $ K_{\kappa,x}(y):=K_\kappa(x,y)$
and denote by
\[
    m(K_{\kappa,x},s)
    =
    w\bigl(\{y\in\mathbb S^n:
    |K_{\kappa,x}(y)|>s\}\bigr),
\]
the distribution function of $K_{\kappa,x}$. Let $K_{\kappa,x}^*$ denote the decreasing rearrangement of
$K_{\kappa,x}$:
\begin{align}
    K_{\kappa,x}^*(t)
    =
    \inf\{s>0:m(K_{\kappa,x},s)\leq t\}.\nonumber
\end{align}
Combining the expansion above with the small-geodesic-ball volume estimate
\[
    w(B_g(x,r))
    =
    \frac{\omega_{n-1}}{n\omega_n}r^n(1+O(r^2))
\]
gives
\begin{align}
    \sup_{x\in\mathbb S^n}
     m(K_{\kappa,x},s)
    \leq\frac{\omega_{n-1}\omega_n}{(2\pi)^nn}
    s^{-2}
    (1+O(s^{-\delta'})),\qquad s\rightarrow+\infty,\nonumber
\end{align}
for some $\delta'>0$.
Set
\begin{align}
     \Phi(s)
    :=
    \frac{\omega_{n-1}\omega_n}{(2\pi)^nn}
    s^{-2}
    (1+O(s^{-\delta'})).\nonumber
\end{align}
For any $t>\Phi(s)$ and $x\in\S^n$, from
\[
    m(K_{\kappa,x},s)\leq \Phi(s)<t
\]
we obtain $ K_{\kappa,x}^*(t)\leq s$. Thus,
\[
    k_1^*(t)
    :=
    \sup_{x\in\mathbb S^n}K_{\kappa,x}^*(t)
    \leq s
\]
whenever $t>\Phi(s)$, and hence
\[
    m(k_1^*,s)\leq \Phi(s).
\]
This is precisely the first distribution estimate required in \cite[Theorem 1]{Fon2}. Since $T_\kappa$ is self-adjoint, the second required distribution estimate follows in the same way:
\begin{align}
    m(k_2^*,s)\leq \Phi(s).\nonumber
\end{align}

Consequently, applying \cite[Theorem 1]{Fon2} to $T_\kappa$ gives
\[
    \int_{\mathbb S^n}
    \exp\left(
       \frac{n!}{2}
       \frac{|T_\kappa f|^2}{\|f\|_{L^2(\S^n)}^2}
    \right)\,\de w
    \leq C(n,\kappa)
\]
for every nonzero function $f\in L^2(\mathbb S^n)$ (here we also use the identity $\omega_{n-1}\omega_{n}=\frac{2^{n+1}\pi^n}{(n-1)!}$). For every $v\in H^\frac{n}{2}(\S^n)$ satisfying
$\int_{\S^n} v\,\mathrm{d}w = 0$, taking $f=\mathcal B_\kappa^{1/2}v$ yields \eqref{penal1}.

Finally, Young's inequality gives
\[
    nv
    \leq
    \frac{n!}{2}\frac{v^2}{Q_\kappa(v)}
    +
    \frac{n^2}{2n!}Q_\kappa(v).
\]
Therefore,
\[
    \log\int_{\mathbb S^n}e^{nv}\,\de w
    \leq
    \frac{n^2}{2n!}Q_\kappa(v)+\log C(n,\kappa)
    =
    \frac{n}{2(n-1)!}Q_\kappa(v)+\log C(n,\kappa),
\]
which is \eqref{penal2}. This completes the proof.
\end{proof}

The second ingredient is the following mass separation lemma.

\begin{lemma}\label{lem4.3}
Assume $\sigma$ is a Borel probability measure on $\S^n$ satisfying
\begin{align}
     \int_{\S^n}f\de\sigma=0
        \quad \text{for every } f\in \PP_{m,n}.\nonumber
\end{align}
Then there exists a positive constant $R=R(m,n)>0$ such that, for any $0<r<\frac{R}{100}$, one can find $N_*$ points $x_1,\cdots,x_{N_*}\in\S^n$ such that $d(x_i,x_j)\geq \frac{R}{2}$ for $i\neq j$ and $\sigma(B(x_i,r))\geq \delta_r$ for every $i$, where $\delta_r=\delta_r(m,n)>0$.
\end{lemma}
\begin{proof}
Let $f_1,\cdots,f_t$ be an orthonormal basis of $\PP_{m,n}$ in $L^2(\de w)$, and define $F=(f_1,\cdots,f_t):\S^n\rightarrow \R^t$. Set $C_1=\norm{F}_{L^\infty(\S^n)}$ and $C_2=\norm{\nabla F}_{L^\infty(\S^n)}$. By the definition of $N_*$ in \eqref{yu} and a standard compactness argument,
\begin{align}
    C_3:=\min\left\{\Big|\sum_{i=1}^{N_*-1} w_iF(x_i)\Big|: w_i\geq 0,\sum_{i=1}^{N_*-1}w_i=1,x_i\in \S^n,1\leq i\leq N_*-1\right\}\nonumber
\end{align}
is a positive number. Set $R=\min\{\frac{\pi}{16},\frac{C_3}{16C_2}\}$ and $\eta=\frac{C_3}{16C_1}$. We claim that for any $N_*-1$ points $x_1,\cdots,x_{N_*-1}\in\S^n$,
\begin{align}
    \sigma\left(\S^n\backslash\bigcup_{i=1}^{N_*-1}B(x_i,R)\right)\geq \eta.\nonumber
\end{align}
Indeed, if not, we take mutually disjoint Borel sets $A_i\subset B(x_i,R)$ such that
\begin{align}
    \bigcup_{i=1}^{N_*-1}A_i=\bigcup_{i=1}^{N_*-1}B(x_i,R)\nonumber
\end{align}
and define a discrete probability measure
\begin{align}
    \sigma_0:=\sum_{i=1}^{N_*-1}\sigma(A_i)\delta_{x_i}+\Big(1-\sum_{i=1}^{N_*-1}\sigma(A_i)\Big)\delta_{x_1}.\nonumber
\end{align}
Since $\int_{\S^n}F\de \sigma=0$, we have
\begin{align}
    \Big|\int_{\S^n}F\de\sigma_0\Big|\leq&\; \sum_{i=1}^{N_*-1}\Big|\int_{A_i}F\de\sigma-\sigma(A_i)F(x_i)\Big|+\Big(1-\sum_{i=1}^{N_*-1}\sigma(A_i)\Big)|F(x_1)|+\Big|\int_{\S^n\backslash\bigcup A_i}F\de\sigma\Big|\nonumber\\
    \leq&\;\sum_{i=1}^{N_*-1}C_2R\sigma(A_i)+2C_1\eta\leq \frac{C_3}{2},\nonumber
\end{align}
which contradicts the definition of $C_3$ above.

To finish the proof, choose finitely many points $y_1,\cdots,y_s$ such that $\S^n\subset\bigcup_i B(y_i,r)$. We may assume that $s\leq Cr^{-n}$ for some positive constant $C=C(n)$. First choose $x_1\in\{y_i:1\leq i\leq s\}$ such that $\sigma(B(x_1,r))\geq s^{-1}$. By the claim above, we can then find $x_2\in\{y_i:1\leq i\leq s\}$ such that $d(x_1,x_2)\geq R-r$ and $\sigma(B(x_2,r))\geq \eta s^{-1}$. Iterating this procedure, we obtain $N_*$ points $x_1,\cdots,x_{N_*}\in\S^n$ such that $d(x_i,x_j)\geq R-r\geq \frac{R}{2}$ for $i\neq j$ and $\sigma(B(x_i,r))\geq\eta s^{-1}$. This completes the proof.
\end{proof}
\begin{remark}\label{rmk4.4}
    The preceding proof also shows that the same conclusion holds under the relaxed condition
    \begin{align}
        \Big|\int_{\S^n}f\de\sigma\Big|\leq \varepsilon_0\norm{f}_{L^\infty}
        \quad \text{for every } f\in \PP_{m,n},\nonumber
    \end{align}
    where $\varepsilon_0=\varepsilon_0(m,n)$ is a small positive number. This observation is useful in the proof of Theorem \ref{thm6.2}.
\end{remark}

The third ingredient provides crucial estimates for suitable orthogonal decompositions of $H^\frac{n}{2}(\S^n)$. Specifically, let $x_1,\cdots,x_{N_*}\in\S^n$ be the points obtained in Lemma \ref{lem4.3}, and assume that $0<r<\frac{R}{100}$. Set $U_i=B(x_i,2r)$ and define
\[
    X_i
    :=
    \overline{
       \left\{
          \varphi\in C_c^\infty(U_i):
          \int_{\mathbb S^n}\varphi\,\de w=0
       \right\}
    }^{\,H^{\frac{n}{2}}(\mathbb S^n)},\quad X:=\bigoplus_{i=1}^{N_*}X_i.
\]
For any $v\in H^\frac{n}{2}(\S^n)$ with $\int_{\mathbb S^n}v\,\de w=0$, there is a unique decomposition $v=\sum\limits_{i=1}^{N_*}v_i+h$ such that $v_i\in X_i$ for every $i$ and $h$ is orthogonal to $X$, that is, $B_n(h,f)=0$ for every $f\in X$. We then have the following estimates.
\begin{lemma}\label{lem4.5}
There exist two positive constants $C_1=C_1(n,R)$ and $C_2=C_2(n,r)$ such that
\begin{align}\label{penal3}
    |B_n(v_i,v_j)|\leq C_1r^n\norm{v_i}_{L^2(\S^n)}\norm{v_j}_{L^2(\S^n)}
\end{align}
for every $i\neq j$, and
\begin{align}\label{penal4}
    \norm{h}_{L^\infty(B(x_i,r))}\leq C_2E_n(h)^\frac{1}{2}
\end{align}
for every $i$.
\end{lemma}
\begin{proof}
We first prove \eqref{penal3}. If \(n\) is even, \(L^{\frac{n}{2}}\) is a local differential operator and preserves
supports. Since the supports of \(v_i\) and \(v_j\) are separated, we have $ B_n(v_i,v_j)=\langle v_i,L^{\frac{n}{2}}v_j\rangle=0$ for every $i\neq j$. Suppose now that \(n\) is odd. Then \(L^{\frac{n}{2}}\) is a classical pseudodifferential operator of order \(n\). Its Schwartz
kernel is smooth away from the diagonal. Hence, on
\(B(x_i,2r)\times B(x_j,2r)\), it is uniformly bounded by a constant $C=C(n,R)$, which further implies that
\begin{align}
    |B_n(v_i,v_j)|
    \leq C(n,R)\|v_i\|_{L^1(\S^n)}\|v_j\|_{L^1(\S^n)}.\nonumber
\end{align}
H\"older's inequality now gives \eqref{penal3}.

We next prove \eqref{penal4}. Choose \(\psi_i\in C_c^\infty(B(x_i,2r))\) with $\int_{\mathbb S^n}\psi_i\,\de w=1$ and put $c_i:=B_n(h,\psi_i)$. For any \(\varphi\in C_c^\infty(B(x_i,2r))\), define
\[
    \varphi_0
    :=
    \varphi-
    \left(\int_{\mathbb S^n}\varphi\,dw\right)\psi_i.
\]
Then \(\varphi_0\in X_i\) and we have
\[
    0
    =
    B_n(h,\varphi_0)
    =
    B_n(h,\varphi)
    -
    c_i\int_{\mathbb S^n}\varphi\,dw.
\]
Thus, $L^{\frac{n}{2}}h=c_i$ weakly in $B(x_i,2r)$. A standard interior estimate for elliptic pseudodifferential operators (see, for instance,
\cite[Chapter~XVIII]{Hor2} or \cite[Theorem~7.2]{Shu}) yields
\begin{align}
     \|h\|_{H^n(B(x_i,r))}
    \leq
    C(n,r)\bigl(\|h\|_{L^2(\S^n)}+|c_i|\bigr).\nonumber
\end{align}
Recall that the Poincar\'e inequality gives $\|h\|_{L^2(\S^n)}^2\leq C(n)E_n(h)$. Moreover, the Cauchy--Schwarz inequality gives
\begin{align}
    |c_i|=|B_n(h,\psi_i)|
    \leq
    E_n(h)^{\frac{1}{2}}E_n(\psi_i)^{\frac{1}{2}}.\nonumber
\end{align}
Choosing $\psi_i$ to depend only on $n$ and $r$, and then applying the embedding \(H^n\hookrightarrow L^\infty\), we obtain \eqref{penal4}.
\end{proof}

We are now ready to prove Theorem \ref{thm4.1}.
\begin{proof}[Proof of Theorem \ref{thm4.1}]
Throughout the proof, $C$ denotes a positive constant depending only on $m$ and $n$, whose value may vary from line to line. Define the measure
\begin{align}
    \de\sigma=\frac{e^{nu}}
    {\int_{\mathbb S^n}e^{nu}\de w}\,\de w\nonumber
\end{align}
on $\S^n$. The assumptions on $u$ imply that $\sigma$ satisfies the hypotheses of Lemma \ref{lem4.3}. We may therefore apply Lemmas \ref{lem4.3} and \ref{lem4.5}, retaining their notation. Choose $0<r=r(m,n)<\frac{R}{100}$, to be fixed sufficiently small below, and write $u=\sum\limits_{i=1}^{N_*}v_i+h$ as in the above decomposition.
Lemma \ref{lem4.3} gives
\begin{align}
    \log\int_{\mathbb S^n}e^{nu}\,\de w
    \leq
    \log\int_{B(x_i,r)}e^{nu}\,\de w-\log\delta_r.\nonumber
\end{align}
As $u=v_i+h$ on $B(x_i,r)$, applying \eqref{penal2} with $\kappa=1$ and \eqref{penal4} yields
\begin{align}
    \log\int_{\mathbb S^n}e^{nu}\,\de w
    &\leq
    \log\int_{\mathbb S^n}e^{nv_i}\,\de w
    +
    n\|h\|_{L^\infty(B(x_i,r))}
    -
    \log\delta_r \notag\\
    &\leq
    \frac{n}{2(n-1)!}
    \bigl(E_n(v_i)-\|v_i\|_{L^2(\S^n)}^2\bigr)
    +
    CE_n(h)^{\frac{1}{2}}
    +
    C.\nonumber
\end{align}
Averaging this estimate over $1\leq i\leq N_*$ gives
\begin{align}
    \log\int_{\mathbb S^n}e^{nu}\,\de w\leq&\; A_{m,n}\sum_{i=1}^{N_*} \bigl(E_n(v_i)-\|v_i\|_{L^2(\S^n)}^2\bigr)+CE_n(h)^{\frac{1}{2}}
    +C\nonumber\\
    =&\;A_{m,n}\Big(E_n(u)-\sum_{i=1}^{N_*}\|v_i\|_{L^2(\S^n)}^2-\sum_{i\neq j}B_n(v_i,v_j)-E_n(h)\Big)+CE_n(h)^{\frac{1}{2}}+C.\nonumber
\end{align}
By choosing $r=r(m,n)$ sufficiently small and combining \eqref{penal3} with H\"older's inequality, we obtain
\begin{align}
    \sum_{i\neq j}|B_n(v_i,v_j)|\leq \sum_{i=1}^{N_*}\|v_i\|_{L^2(\S^n)}^2.\nonumber
\end{align}
Therefore,
\begin{align}
    \log\int_{\mathbb S^n}e^{nu}\,\de w\leq&\; A_{m,n}\Big(E_n(u)-E_n(h)\Big)+CE_n(h)^{\frac{1}{2}}+C\nonumber\\
    =&\;A_{m,n}E_n(u)-A_{m,n}E_n(h)+CE_n(h)^\frac{1}{2}+C\nonumber\\
    \leq&\; A_{m,n}E_n(u)+C.\nonumber
\end{align}
The proof is complete.
\end{proof}

\subsection{The best remainder}
In this subsection, we set
\begin{align}
     C_{\mathrm{best}}(m,n)
    :=
    \sup_{u}
    \left\{
       \log\int_{\mathbb S^n}e^{nu}\,\de w
       -
       A_{m,n}E_n(u)
    \right\},\nonumber
\end{align}
where the supremum is taken over all $u\in H^\frac{n}{2}(\S^n)$ satisfying $\int_{\S^n} u\,\mathrm{d}w = 0$ and \eqref{or2}. This definition is equivalent to the one used in Theorem \ref{thm1.5}.

Our main result is the following, which is slightly stronger than Theorem \ref{thm1.5}.
\begin{theorem}\label{thm4.6}
    Assume $m\geq 1$ and $n\geq2$. If
    \begin{align}\label{pos1}
        [(m+1)(m+n)]^\frac{n}{2}< N_*n!,
    \end{align}
    then $C_{\mathrm{best}}(m,n)>0$. In particular, \eqref{pos1} holds in the following four cases: $m=1,n\geq 9$; $m=2,n\geq 12$; $m=3,n\geq 18$; $m\geq 4,n\geq 18(m+1)$.

    Moreover, $C_{\mathrm{best}}(m,n)>0$ when $m=1$ and $5\leq n\leq 8$.
\end{theorem}
The proof of Theorem \ref{thm4.6} has two parts. In the first, we derive \eqref{pos1} by constructing test functions that are suitable perturbations of constant functions. In the second, we treat the case $m=1$ and $5\leq n\leq 8$ by constructing explicit test functions.
\begin{proof}[Proof of the first part]
    Let $Y$ be a nonzero function in $\mathcal{H}_{m+1}$. For $|t|$ sufficiently small, define the smooth function $u_t=\frac{1}{n}\log (1+tY)-\frac{1}{n}\int_{\S^n}\log (1+tY)\de w$. It is straightforward to verify that $\int_{\S^n} u_t\,\mathrm{d}w = 0$ and that $u_t$ satisfies \eqref{or2}. Thus, $u_t$ is an admissible test function.

    A Taylor expansion gives
    \begin{align}
        \log\int_{\mathbb S^n}e^{nu_t}\,\de w=\frac{t^2}{2}\norm*{Y}_{L^2(\S^n)}^2+O(t^3).\nonumber
    \end{align}
    Moreover, since $u_t=\frac{t}{n}Y+O(t^2)$ and $Y\in\mathcal{H}_{m+1}$, we have
    \begin{align}
        E_n(u_t)=\frac{t^2}{n^2}E_n(Y)+O(t^3)=\frac{t^2}{n^2}\lambda_{m+1}^\frac{n}{2}\norm*{Y}_{L^2(\S^n)}^2+O(t^3).\nonumber
    \end{align}
    Combining these two estimates and recalling that $\lambda_{m+1}=(m+1)(m+n)$, we obtain
    \begin{align}
        \log\int_{\S^n}e^{nu_t}\,\de w-A_{m,n}E_n(u_t)=&\;\Big(\frac{1}{2}-\frac{A_{m,n}}{n^2}[(m+1)(m+n)]^\frac{n}{2}\Big)t^2\norm*{Y}_{L^2(\S^n)}^2+O(t^3)\nonumber\\
        =&\;\Big(N_*n!-[(m+1)(m+n)]^\frac{n}{2}\Big)\frac{A_{m,n}}{n^2}t^2\norm*{Y}_{L^2(\S^n)}^2+O(t^3).\nonumber
    \end{align}
    Therefore, \eqref{pos1} implies that $C_{\mathrm{best}}(m,n)>0$.

    Finally, since $\Theta(1,0,n)=2$, $\Theta(2,0,n)=n+2$, and $\Theta(3,0,n)=2n+2$, substituting these explicit values into \eqref{pos1} and performing elementary calculations shows that \eqref{pos1} holds in the following three cases: $m=1,n\geq 9$; $m=2,n\geq 12$; and $m=3,n\geq 18$. For $m\geq 4$ and $n\geq 18(m+1)$, we have
    \begin{align}
        n!>\frac{n^n}{e^n}>\frac{n^n}{3^n}\geq [2(m+1)n]^\frac{n}{2}\geq [(m+1)(n+m)]^\frac{n}{2},\nonumber
    \end{align}
    which also implies that \eqref{pos1} holds.
\end{proof}

\begin{proof}[Proof of the second part]
    Put $z=x_{n+1}$ and $\alpha=\frac{n-1}{2}$, and let \(C_k^\alpha\) be the standard Gegenbauer polynomial, normalized by
    \begin{align}
        C_0^\alpha=1,
    \qquad
    C_1^\alpha=2\alpha z,\nonumber
    \end{align}
    and
    \begin{align}
        (k+1)C_{k+1}^\alpha
    =
    2(k+\alpha)zC_k^\alpha
    -
    (k+2\alpha-1)C_{k-1}^\alpha.\nonumber
    \end{align}
    We define
    \begin{align}
    u_5
    ={}&
    \frac{1789}{10000}C_2^2
    +
    \frac{41}{2000}C_4^2
    +
    \frac{3}{625}C_6^2
    +
    \frac{1}{625}C_8^2
    +
    \frac{3}{5000}C_{10}^2
    +
    \frac{1}{5000}C_{12}^2\nonumber\\
    ={}&
    \frac{6656}{625}z^{12}
    -
    \frac{12672}{625}z^{10}
    +
    \frac{9504}{625}z^8
    -
    \frac{2912}{625}z^6+
    \frac{169}{125}z^4
    +
    \frac{3681}{2500}z^2
    -
    \frac{3097}{10000},\nonumber\\
    u_6
    ={}&
    \frac65z^4+\frac85z^2-\frac27
    =
    \frac{248}{1925}C_2^{5/2}
    +
    \frac{16}{1925}C_4^{5/2},\nonumber\\
    u_7
    ={}&
    3\left(z^2-\frac18\right),
    \qquad
    u_8
    =
    4\left(z^2-\frac19\right).\nonumber
\end{align}
We refer to \cite{Sze} for background on Gegenbauer polynomials. Since \(C_k^\alpha(z)\in\mathcal H_k\), it is straightforward to verify that $u_n$ satisfies $\int_{\S^n} u_n\,\mathrm{d}w = 0$ and \eqref{or2} for every $5\leq n\leq 8$. Moreover, direct computations give
\begin{align}
    A_{1,5}E_5(u_5)<\frac{225}{32},\quad A_{1,6}E_6(u_6)=\frac{17504}{3003},\quad A_{1,7}E_7(u_7)=\frac{196}{25},\quad A_{1,8}E_8(u_8)=\frac{4608}{385}.\nonumber
\end{align}
To obtain lower bounds for $\int_{\mathbb S^n}e^{nu_n}\,\de w$, we use the elementary estimate
\begin{align}
     T_N(x):=\sum_{j=0}^N\frac{x^j}{j!}\leq e^x,
    \qquad \forall\;x\in\mathbb R,
    \quad \forall\;N\ \text{odd}.\nonumber
\end{align}
Set $I_{n,N}(u)=\int_{\mathbb S^n}T_N(nu)\,\de w$. Direct computations yield
\begin{align}
    I_{5,19}(u_5)\geq 1200,\quad I_{6,13}(u_6)\geq 400,\quad I_{7,13}(u_7)\geq3000,\quad I_{8,13}(u_8)\geq200000.\nonumber
\end{align}
Since
\begin{align}
    e^\frac{225}{32}<1200,\quad e^{\frac{17504}{3003}}<400,\quad e^\frac{196}{25}<3000,\quad e^\frac{4608}{385}<200000,\nonumber
\end{align}
we conclude that $C_{\mathrm{best}}(1,n)>0$ for $5\leq n\leq 8$.
\end{proof}

\section{Quantitative Sobolev remainder estimates}\label{sec5}
Throughout this section, we fix $m,n\in\N_+$ and $1\leq p<n$, and set
$p^*=\frac{np}{n-p}$ and $\theta=\frac{p}{p^*}=\frac{n-p}{n}$. For
notational simplicity, we write
\[
    \Theta=\Theta(m,\theta,n).
\]
Our main result is the following, which is stronger than Theorem~\ref{thm1.6}.
\begin{theorem}\label{thm5.1}
Assume that the \textbf{sharp} Sobolev inequality \eqref{Sob3} holds on
$\S^n$. Then there exists a positive constant $B=B(m,n,p)$ such that, for
every $0<\varepsilon<1$, every linear functional
$T:\PP_{m,n}\rightarrow\R$, and every $u\in W^{1,p}(\S^n)$ satisfying
\begin{align}\label{eq5.moment-condition}
    \left|\int_{\S^n}f|u|^{p^*}\de\mu\right|
    \leq |T(f)|\norm*{u}_{L^p(\S^n)}^{p^*}
    \quad \text{for all } f\in\PP_{m,n},
\end{align}
one has
\begin{align}\label{eq5.main}
    \norm*{u}_{L^{p^*}(\S^n)}^p
    \leq
    \left(
        \frac{S_{n,p}^p}{\Theta(m,\frac{n-p}{n},n)}+\varepsilon
    \right)
    \norm*{\nabla u}_{L^p(\S^n)}^p
    +B\left(\varepsilon^{-\gamma(p)}+\norm{T}_\infty^\theta\right)
    \norm*{u}_{L^p(\S^n)}^p,
\end{align}
where $\gamma(p)$ is defined in \eqref{newcst}, and
\begin{align}\label{eq5.beta}
    \norm{T}_\infty
    =\sup_{\substack{f\in\PP_{m,n}\\
    \norm{f}_{L^\infty(\S^n)}\leq1}}|T(f)|.
\end{align}
\end{theorem}

The proof of Theorem~\ref{thm5.1} requires several quantitative estimates.
Let $\de w=\frac{\de\mu}{\omega_n}$ denote the normalized surface measure on
$\S^n$.
We set
\begin{align}
    D=\dim(\P_{m,n}|_{\S^n})
    =\binom{n+m}{n}+\binom{n+m-1}{n},
    \qquad d=D-1.
\end{align}
Choose an orthonormal basis $f_1,\ldots,f_d$ of $\PP_{m,n}$ in
$L^2(\de w)$ and define $F=(f_1,\ldots,f_d)$. The rotational invariance of
$\PP_{m,n}$ shows that $|F(x)|^2$ is constant on $\S^n$. Consequently,
\begin{align}\label{eq5.frame-identities}
    |F(x)|^2=d,\qquad
    \int_{\S^n}F\,\de w=0,\qquad
    \int_{\S^n}FF^{\mathsf T}\,\de w=\operatorname{Id}.
\end{align}
For a probability measure $\sigma$ on $\S^n$, define
\begin{align}
    \tau(\sigma)
    =\sup_{\substack{f\in\PP_{m,n}\\
    \norm{f}_{L^\infty(\S^n)}\leq1}}
    \left|\int_{\S^n}f\,\de\sigma\right|.
\end{align}
If $e=\int_{\S^n}F\,\de\sigma$, then we have
\begin{align}\label{eq5.e-bound}
    |e|\leq\sqrt d\,\tau(\sigma).
\end{align}
The first lemma gives a lower bound for $\norm{\cdot}_\theta$ on general
discrete measures.
\begin{lemma}\label{lem5.2}
For every finitely supported probability measure
$\nu=\sum
\limits_i\nu_i\delta_{x_i}$ on $\S^n$,
\begin{align}\label{eq5.moment-correction}
    \sum_i\nu_i^\theta
    \geq \Theta-D^{1-\theta}d^\theta\tau(\nu)^\theta.
\end{align}
\end{lemma}
\begin{proof}
Let $e=\int_{\S^n}F\,\de\nu$. If $e=0$, then \eqref{eq5.moment-correction}
follows directly from the definition of $\Theta$ in \eqref{Cub1}. Suppose
that $e\neq0$. We correct $\nu$ by adding a suitable measure. Consider the
probability measure
\[
    \de\rho(x)
    =\left(1-\frac{e\cdot F(x)}{\sqrt d\,|e|}\right)\de w(x).
\]
The density is nonnegative by \eqref{eq5.frame-identities}, and the same
identities give
$\int_{\S^n}F\,\de\rho=-\frac{e}{\sqrt d\,|e|}$. By Carath\'eodory's
theorem, there is a probability measure $\rho_a$ supported on at most
$d+1=D$ points and satisfying
$\int_{\S^n}F\,\de\rho_a=\int_{\S^n}F\,\de\rho$.

Set $t=\sqrt d\,|e|$ and
$\widehat\nu=(\nu+t\rho_a)/(1+t)$. Then $\widehat\nu\in\M^c_{m,n}$. Since
$(a+b)^\theta\leq a^\theta+b^\theta$ for $a,b\geq0$, we have
\[
    \Theta
    \leq\norm{\widehat\nu}_\theta
    \leq\frac{\norm{\nu}_\theta+t^\theta\norm{\rho_a}_\theta}
    {(1+t)^\theta}
    \leq\norm{\nu}_\theta+D^{1-\theta}t^\theta.
\]
By \eqref{eq5.e-bound}, $t\leq d\tau(\nu)$, which proves the lemma.
\end{proof}

Let $q=\theta^{-1}=\frac{p^*}{p}$. Given finitely many continuous functions
$a_j\geq0$ with $\sum\limits_j a_j=1$, define (the specific choice of the
$a_j$ will be made in the proof of Theorem~\ref{thm5.1})
\[
    \Phi(\sigma)=\sum_j\left(\int_{\S^n}a_j^q\,\de\sigma\right)^\theta.
\]
The following lemma provides a finite-atom reduction principle.
\begin{lemma}\label{lem5.3}
For every probability measure $\sigma$ on $\S^n$, there is a probability
measure $\nu=\sum\limits_{i=1}^K w_i\delta_{x_i}$ such that
\begin{align}\label{eq5.atomic-reduction}
    K\leq D,\qquad w_i>0,\qquad
    \int_{\S^n}F\,\de\nu=\int_{\S^n}F\,\de\sigma,
    \qquad \Phi(\nu)\leq\Phi(\sigma).
\end{align}
\end{lemma}

\begin{proof}
Apply Carath\'eodory's theorem to the continuous map
$(F,a_1^q,\ldots,a_J^q)$ to replace $\sigma$ by a finitely supported
probability measure $\nu_0$ preserving the integrals of all its components.
Suppose that the resulting measure has $K>D=d+1$ positive weights.
There is then a nonzero vector $(h_i)_{i=1}^K$ such that
\[
    \sum_i h_i=0,\qquad \sum_i h_iF(x_i)=0.
\]
The set of parameters $s$ for which $w_i+sh_i\geq0$ for every $i$ is a closed
interval $[s_-,s_+]$ with $s_-<0<s_+$, and at each endpoint at least one
positive weight vanishes. Moreover, the function
\[
    s\longmapsto
    \sum_j\left(\sum_i(w_i+sh_i)a_j(x_i)^q\right)^\theta
\]
is concave. Hence its value at one of the two endpoints does not exceed its
value at $s=0$. Choosing that endpoint preserves the mass and the $F$-moment,
does not increase $\Phi$, and decreases the number of positive weights.
Iteration gives \eqref{eq5.atomic-reduction}.
\end{proof}

Our third lemma constructs a suitable $\ell^p$-partition of unity.
\begin{lemma}\label{lem5.4}
There exist positive constants $c,C,r_1$, depending only on $n$ and $p$, such
that for every $0<r<r_1$ there are finitely many nonnegative $C^1$ functions
$\eta_j$ supported in geodesic balls of radius $r$ for which, with
$a_j=\eta_j^p$,
\begin{align}\label{eq5.partition}
    \sum_j\eta_j^p=1,\qquad
    \sum_j|\nabla\eta_j|^p\leq Cr^{-p},
\end{align}
and
\begin{align}\label{eq5.resolution}
    \sum_j|a_j(x)-a_j(y)|^2
    \geq c\min\left\{\frac{|x-y|^2}{r^2},1\right\}
    \quad \text{for all }x,y\in\S^n.
\end{align}
\end{lemma}
\begin{proof}
Take a maximal $\frac{r}{8}$-separated set $\{z_\ell\}$ and choose nonnegative smooth
cutoff functions $\psi_\ell$ which are equal to one on
$B(z_\ell,\frac{r}{2})$, are
supported in $B(z_\ell,r)$, and satisfy
$|\nabla\psi_\ell|\leq \frac{C}{r}$. These balls have overlap bounded by a
constant $M=M(n)$. Set
\[
    \eta^0_\ell
    =\frac{\psi_\ell}{(\sum\limits_k\psi_k^p)^{1/p}},
    \qquad \lambda_\ell=(\eta^0_\ell)^p.
\]
Then $\sum\limits_\ell\lambda_\ell=1$,
$\sum\limits_\ell|\nabla\eta^0_\ell|^p\leq Cr^{-p}$, and for every $x\in\S^n$
there is an index $\ell$ such that $\lambda_\ell(x)\geq M^{-1}$. To obtain
\eqref{eq5.resolution}, we modify the functions $\eta_\ell^0$ as follows.

Set $N=n+1$ and choose a smooth function $\chi$ such that
$\chi(t)=t$ for $|t|\leq1/2$, $|\chi|\leq3/4$, and $|\chi'|\leq C$. For
$1\leq k\leq N$ and each sign, let
\[
    h_{\ell,k,\pm}(x)
    =\frac{1\pm\chi((x_k-z_{\ell,k})/(2r))}{2N},
    \qquad
    \eta_{\ell,k,\pm}=\eta^0_\ell h_{\ell,k,\pm}^{1/p}.
\]
Since $\sum\limits_{k,\pm}h_{\ell,k,\pm}=1$,
$h_{\ell,k,\pm}\geq(8N)^{-1}$, and
$|\nabla h_{\ell,k,\pm}|\leq \frac{C}{r}$, the functions
$\eta_{\ell,k,\pm}$ satisfy \eqref{eq5.partition}. Denote their $p$-th powers
by $a_{\ell,k,\pm}=\lambda_\ell h_{\ell,k,\pm}$.

It remains to establish \eqref{eq5.resolution}. Fix $x,y\in\S^n$ and choose
$\ell$ so that $\lambda_\ell(x)\geq M^{-1}$. Let $\Delta a$ be the vector
consisting of the $2N$ differences
$a_{\ell,k,\pm}(x)-a_{\ell,k,\pm}(y)$. If $\lambda_\ell(y)=0$, then the sum
of these differences is $\lambda_\ell(x)$, and hence
$\norm{\Delta a}_{\ell^2}\geq(M\sqrt{2N})^{-1}$. If
$\lambda_\ell(y)>0$, then $x,y\in B(z_\ell,r)$, so $\chi$ is linear at all
the arguments involved. With $A=\lambda_\ell(x)-\lambda_\ell(y)$, one has
\[
    A=\sum_{k,\pm}\Delta a_{k,\pm},
    \qquad
    2Nr(\Delta a_{k,+}-\Delta a_{k,-})
    =\lambda_\ell(x)(x_k-y_k)+A(y_k-z_{\ell,k}).
\]
Since $|y-z_\ell|\leq r$, H\"older's inequality gives
$\lambda_\ell(x)|x-y|\leq C(n)r\norm{\Delta a}_{\ell^2}$. Combining the two
cases and adding the nonnegative contributions corresponding to the other
indices proves \eqref{eq5.resolution}. Relabeling the triples
$(\ell,k,\pm)$ as a single index concludes the proof.
\end{proof}
The fourth lemma provides an elementary but useful inequality.
\begin{lemma}\label{lem5.5}
Let $0\leq X\leq1$ be a random variable on a probability space and let $q>1$.
Then
\begin{align}\label{eq5.convexity}
    \norm{X}_{L^q}-\mathbb E X
    \geq c_q\operatorname{Var}(X),
    \qquad c_q=\frac12\min\{q-1,1\}.
\end{align}
\end{lemma}

\begin{proof}
If $q\geq2$, then $\norm{X}_{L^q}\geq\norm{X}_{L^2}$ and
\[
    \norm{X}_{L^2}-\mathbb E X
    =\frac{\operatorname{Var}(X)}{\norm{X}_{L^2}+\mathbb E X}
    \geq\frac12\operatorname{Var}(X).
\]
If $1<q<2$, the second derivative of $t\mapsto t^q$ on $(0,1]$ is at least
$q(q-1)$, so Taylor's theorem gives
\[
    \mathbb E X^q
    \geq(\mathbb E X)^q+\frac{q(q-1)}2\operatorname{Var}(X).
\]
The elementary inequality
$(a+b)^\frac{1}{q}\geq a^\frac{1}{q}+\frac{b}{q}$ for
$a,b\geq0$ and $a+b\leq1$ now yields \eqref{eq5.convexity}.
\end{proof}
The next lemma gives a useful error estimate comparing a measure with
small support to a suitably chosen point mass.
\begin{lemma}\label{lem5.6}
Let $\pi$ be a probability measure supported in a spherical region of small
diameter $d_0=d_0(m,n)$. Define
\begin{align}
    b=\int_{\S^n}x\,\de\pi(x),\qquad
    z=\frac{b}{|b|},\qquad
    v=1-|b|^2
    =\frac12\iint_{\S^n\times\S^n}|x-y|^2\,\de\pi(x)\de\pi(y).
\end{align}
Then, for every $f\in\P_{m,n}$,
\begin{align}
    \left|\int_{\S^n}f\,\de\pi-f(z)\right|
    \leq C(m,n)\norm{f}_{L^\infty(\S^n)}v.
\end{align}
\end{lemma}
\begin{proof}
Let $C=C(m,n)$ be a positive constant that may vary from line to line. Since
$\P_{m,n}$ is finite-dimensional, we have
\[
    \norm{f}_{C^2(\S^n)}
    \leq C\norm{f}_{L^\infty(\S^n)}.
\]
Taking $d_0$ sufficiently small and using the Taylor expansion of $f$ at $z$
together with the definition of $z$, we obtain
\begin{align}
    \left|\int_{\S^n}f\,\de\pi-f(z)\right|
    \leq&\; C\norm{f}_{L^\infty(\S^n)}\Big(\big|\int_{\S^n}(x-z)\de\pi(x)\big|+\int_{\S^n}|x-z|^2\,\de\pi(x)\Big)\nonumber\\
    \leq&\;C\norm{f}_{L^\infty(\S^n)}(1-|b|)\leq C\norm{f}_{L^\infty(\S^n)}v.\nonumber
\end{align}
The proof is complete.
\end{proof}

For the partition in Lemma~\ref{lem5.4}, define
\begin{align}\label{eq5.Phi-r}
    \Phi_r(\sigma)
    =\sum_j\left(\int_{\S^n}\eta_j^{p^*}\,\de\sigma\right)^\theta
    =\sum_j\left(\int_{\S^n}a_j^q\,\de\sigma\right)^{1/q}.
\end{align}
Combining the five lemmas above, we obtain the following lower bound for
$\Phi_r(\sigma)$, which is crucial for proving \eqref{eq5.main}.
\begin{proposition}\label{prop5.7}
There exist $r_0=r_0(m,n,p)>0$ and $C=C(m,n,p)>0$ such that, for every
$0<r<r_0$ and every probability measure $\sigma$ on $\S^n$,
\begin{align}\label{eq5.localization}
    \Phi_r(\sigma)
    \geq\Theta-Cr^{\frac{2\theta}{1-\theta}}-C\tau(\sigma)^\theta.
\end{align}
\end{proposition}

\begin{proof}
By Lemma~\ref{lem5.3}, there is a probability measure
$\nu=\sum\limits_{i=1}^K w_i\delta_{x_i}$ with $K\leq D$,
$\tau(\nu)=\tau(\sigma)$, and $\Phi_r(\nu)\leq\Phi_r(\sigma)$. Thus it suffices to establish a lower bound for $\Phi_r(\nu)$. Join two
atoms $x_i,x_j$ if $d(x_i,x_j)\leq 2r$, and denote the connected
components by $I_1,\ldots,I_L$. Each component has diameter at most
$2(D-1)r$. We choose $r_0<\min\{r_1,\frac{d_0}{2(D-1)}\}$ so that every component lies in
the range of Lemma~\ref{lem5.6}. Since each $a_j$ is supported in a
ball of radius $r$, it can be nonzero at atoms belonging to at most one
component. For each component, set
\[
    A_\ell=\sum_{i\in I_\ell}w_i,\qquad
    \pi_\ell=A_\ell^{-1}\sum_{i\in I_\ell}w_i\delta_{x_i},\qquad
    b_\ell=\int_{\S^n}x\,\de\pi_\ell(x),
\]
and let $z_\ell=\frac{b_\ell}{|b_\ell|}$ and $v_\ell=1-|b_\ell|^2$. Define
\begin{align}
    \overline\nu=\sum_{\ell=1}^L A_\ell\delta_{z_\ell},\qquad
    H=\sum_{\ell=1}^L A_\ell^\theta,\qquad
    V=\sum_{\ell=1}^L A_\ell v_\ell.
\end{align}

Applying Lemma~\ref{lem5.5} to each $a_j$ on each component and using
$\sum\limits_j\int_{\S^n} a_j\,\de\pi_\ell=1$, we obtain
\begin{align}\label{eq5.convexity-gain}
    \sum_{j\,:\,a_j|_{I_\ell}\not\equiv0}
    \left(\int_{\S^n}a_j^q\,\de\nu\right)^\theta=A_\ell^\theta\sum_{j}\norm{a_j}_{L^q(\pi_\ell)}
    \geq
    A_\ell^\theta
    \left(1+c_q\sum_j\operatorname{Var}_{\pi_\ell}(a_j)\right).
\end{align}
By \eqref{eq5.resolution} and the diameter bound on $I_\ell$,
\begin{align*}
    \sum_j\operatorname{Var}_{\pi_\ell}(a_j)
    &=\frac12\iint\sum_j|a_j(x)-a_j(y)|^2
      \,\de\pi_\ell(x)\de\pi_\ell(y)\\
    &\geq\frac{c}{D^2r^2}\,
      \frac12\iint|x-y|^2\,\de\pi_\ell(x)\de\pi_\ell(y)
    =\frac{c}{D^2r^2}v_\ell.
\end{align*}
Summing \eqref{eq5.convexity-gain} over $\ell$ and using
$A_\ell^\theta\geq A_\ell$ gives
\begin{align}\label{eq5.Phi-gain}
    \Phi_r(\sigma)\geq\Phi_r(\nu)
    \geq H+cr^{-2}\sum_\ell A_\ell^\theta v_\ell
    \geq H+cr^{-2}V.
\end{align}

By Lemma~\ref{lem5.6}, for every $f\in\PP_{m,n}$ with
$\norm{f}_{L^\infty(\S^n)}\leq1$,
\[
    \left|\int_{\S^n}f\,\de\overline\nu
    -\int_{\S^n}f\,\de\nu\right|
    \leq C\sum_\ell A_\ell v_\ell=CV.
\]
Thus $\tau(\overline\nu)\leq\tau(\sigma)+CV$. Applying
Lemma~\ref{lem5.2} to $\overline\nu$, we find
\begin{align}\label{eq5.H-lower}
    H\geq\norm{\overline\nu}_\theta
    \geq\Theta-C(\tau(\sigma)+V)^\theta
    \geq\Theta-C\tau(\sigma)^\theta-CV^\theta.
\end{align}
Combining \eqref{eq5.Phi-gain} and \eqref{eq5.H-lower} yields
\begin{align}\label{eq5.localization-strong}
    \Phi_r(\sigma)
    \geq\Theta-C\tau(\sigma)^\theta-CV^\theta+cr^{-2}V.
\end{align}
For positive constants $C_1,c_1$, we have
\[
    \sup_{V\geq0}\{C_1V^\theta-c_1r^{-2}V\}
    =(1-\theta)C_1
    \left(\frac{\theta C_1}{c_1}\right)^{\frac{\theta}{1-\theta}}
    r^{\frac{2\theta}{1-\theta}}.
\]
Optimizing the last two terms in \eqref{eq5.localization-strong} proves
\eqref{eq5.localization}.
\end{proof}
We now prove Theorem~\ref{thm5.1}.
\begin{proof}[Proof of Theorem~\ref{thm5.1}]
Let $B=B(n,p)$ be a fixed constant for which \eqref{Sob3} holds on
$\S^n$, and take the partition from Lemma~\ref{lem5.4}. When
$1\leq p\leq2$, we apply the elementary estimate
\[
    |a+b|^p\leq|a|^p+p|a|^{p-2}a\cdot b+C_p|b|^p
\]
with $a=\eta_j\nabla u$ and $b=u\nabla\eta_j$. The first-order terms cancel
after summation because
\begin{align}\label{eq5.cross-cancellation}
    p|\nabla u|^{p-2}u\nabla u\cdot
    \sum_j\eta_j^{p-1}\nabla\eta_j
    =\frac1p p|\nabla u|^{p-2}u\nabla u\cdot\nabla\sum_j\eta_j^p=0.
\end{align}
Consequently, we have
\begin{align}\label{eq5.energy-low-p}
    \sum_j\norm*{\nabla(\eta_j u)}_{L^p(\S^n)}^p
    \leq\norm*{\nabla u}_{L^p(\S^n)}^p
    +Cr^{-p}\norm*{u}_{L^p(\S^n)}^p.
\end{align}

For $p>2$, we combine the vector inequality
\[
    |a+b|^p
    \leq |a|^p+p|a|^{p-2}a\cdot b
    +C_p\bigl(|a|^{p-2}|b|^2+|b|^p\bigr)
\]
with Young's inequality, valid for $0<\delta\leq1$,
\[
    C_p|a|^{p-2}|b|^2
    \leq\delta|a|^p+C_p\delta^{-\frac{p-2}{2}}|b|^p.
\]
Using the cancellation in \eqref{eq5.cross-cancellation}, we similarly obtain
\begin{align}\label{eq5.energy-high-p}
    \sum_j\norm*{\nabla(\eta_j u)}_{L^p(\S^n)}^p
    \leq(1+\delta)\norm*{\nabla u}_{L^p(\S^n)}^p
    +C\delta^{1-\frac p2}r^{-p}\norm*{u}_{L^p(\S^n)}^p.
\end{align}

The assertion is immediate for $u=0$, so assume that $u\neq0$ and set
\[
    U=\norm*{u}_{L^{p^*}(\S^n)}^p,\qquad
    E=\norm*{\nabla u}_{L^p(\S^n)}^p,\qquad
    L=\norm*{u}_{L^p(\S^n)}^p,\qquad
    \de\sigma=\frac{|u|^{p^*}\de\mu}
    {\norm*{u}_{L^{p^*}(\S^n)}^{p^*}}.
\]
It follows from \eqref{eq5.moment-condition} that
\begin{align}\label{eq5.tau-T}
    \tau(\sigma)^\theta U\leq\norm{T}_\infty^\theta L.
\end{align}
Applying \eqref{Sob3} to each $\eta_j u$ and summing, we obtain
\begin{align}\label{eq5.localized-Sobolev}
    \Phi_r(\sigma)U
    \leq S_{n,p}^p
    \sum_j\norm*{\nabla(\eta_j u)}_{L^p(\S^n)}^p+BL.
\end{align}
Combining Proposition~\ref{prop5.7} with \eqref{eq5.tau-T} and
\eqref{eq5.localized-Sobolev}, and setting
\begin{align}
    \alpha=\frac{2\theta}{1-\theta}=\frac{2(n-p)}p,
\end{align}
we obtain
\begin{align}
    (\Theta-Cr^\alpha)U
    \leq S_{n,p}^p
    \sum_j\norm*{\nabla(\eta_j u)}_{L^p(\S^n)}^p
    +(B+C\norm{T}_\infty^\theta)L.
\end{align}
When $1\leq p\leq 2$, take
$r=c\varepsilon^\frac{1}{\alpha}$ for a sufficiently small positive constant
$c=c(m,n,p)$ and apply \eqref{eq5.energy-low-p} to obtain
\eqref{eq5.main}. When $p>2$, take
$r=c\varepsilon^\frac{1}{\alpha}$ as above and apply
\eqref{eq5.energy-high-p} with $\delta=c\varepsilon$ to obtain the same
conclusion.
\end{proof}
We end this section with the proof of Corollary~\ref{cor1.8}.
\begin{proof}[Proof of Corollary \ref{cor1.8}]
    For simplicity, let $C$ denote a constant depending only on $n,p,m$;
    its value may change from line to line.

    For any $0<\varepsilon<1$, Theorem~\ref{thm1.6} gives
    \begin{align}\label{bb1}
     \norm{u}_{L^{p^*}(\S^n)}^p
    \leq
    \Bigg(
        \frac{S_{n,p}^p}{\Theta\left(m,\frac{n-p}{n},n\right)}+\varepsilon
    \Bigg)
    \norm{\nabla u}_{L^p(\S^n)}^p
    +C\varepsilon^{-\gamma(p)}\norm{u}_{L^{p}(\S^n)}^p.
\end{align}
    We consider the case
    $\norm{\nabla u}_{L^p(\S^n)}>\norm{u}_{L^{p}(\S^n)}$. Choose
    \begin{align}
        \varepsilon=\left(\frac{\norm{u}_{L^{p}(\S^n)}}{\norm{\nabla u}_{L^p(\S^n)}}\right)^\frac{p}{1+\gamma(p)}
        =\left(\frac{\norm{u}_{L^{p}(\S^n)}}{\norm{\nabla u}_{L^p(\S^n)}}\right)^{\Bar{\gamma}(p)},\nonumber
    \end{align}
    It follows from \eqref{bb1} that
    \begin{align}\label{bb2}
     \norm{u}_{L^{p^*}(\S^n)}^p
    \leq
        \frac{S_{n,p}^p}{\Theta\left(m,\frac{n-p}{n},n\right)}
    \norm{\nabla u}_{L^p(\S^n)}^p
    +C \norm{\nabla u}_{L^p(\S^n)}^{p-\Bar{\gamma}(p)}\norm{u}_{L^{p}(\S^n)}^{\Bar{\gamma}(p)}.
\end{align}
    We use the elementary inequality
    \begin{align}\label{bb3}
        (x+y)^q\geq x^q+qx^{q-1}y,\quad x,y\geq0,\quad q\geq1,
    \end{align}
    with $q=1+\gamma(p)$,
    \begin{align}
        x=\Bigg( \frac{S_{n,p}}{\Theta\left(m,\frac{n-p}{n},n\right)^{\frac{1}{p}}} \Bigg)^{\Bar{\gamma}(p)} \norm*{\nabla u}_{L^p(\S^n)}^{\Bar{\gamma}(p)}\nonumber
    \end{align}
    and
    \begin{align}
        y=C\norm*{u}_{L^{p}(\S^n)}^{\Bar{\gamma}(p)},\nonumber
    \end{align}
    we see that the right-hand side of \eqref{bb2} is bounded by $(x+y)^q$. Hence,
    \begin{align}
        \norm{u}_{L^{p^*}(\S^n)}^{\frac{p}{q}} \le x+y,\nonumber
    \end{align}
    which yields \eqref{scaline3}.
\end{proof}

\section{Quantitative Moser--Trudinger--Onofri remainder estimates}
In the first subsection, we prove Theorems \ref{thm1.7} and \ref{thm6.2}, as well as Corollary \ref{cor1.9}. In the second subsection, we generalize \cite[Theorem 1.1]{Chenzhang}.
\subsection{General moment constraints}\label{sec6.1}
The proofs of Theorems \ref{thm1.7} and \ref{thm6.2} rely on the following quantitative ``distribution of mass''
principle for the higher-order Moser--Trudinger--Onofri inequality \eqref{MTO0}.
\begin{lemma}\label{lem6.1}
    Let $(M,g)$ be a smooth closed manifold of dimension $n\geq2$.
    Assume $s\in\N_+$, $1\leq s<n$, and $u\in W^{s,\frac{n}{s}}(M)$ satisfies $\int_M u \de\mu=0$.
    Suppose there exist open subsets $\Omega_1,\cdots,\Omega_k\subset M$ and positive constants $r<1$ and $\delta$ such that $\text{dist}(\Omega_i,\Omega_j)\geq 2r$ whenever $i\neq j$, and
    \begin{align}\label{distMTO1}
        \frac{\int_{\Omega_i}e^{nu}\de\mu}{\int_{M}e^{nu}\de\mu}\geq \delta
    \end{align}
    for all $1\leq i\leq k$. Then for any $0<\varepsilon<1$, we have
    \begin{align}\label{distMTO2}
        \log \int_{M} e^{nu}\,\mathrm{d}\mu
    \leq
    \left(
        \frac{\alpha_{s,n}}{k}
        +\varepsilon
    \right)
    \norm{\Delta^{\frac{s}{2}} u}_{L^{\frac{n}{s}}(M)}^{\frac{n}{s}}
    +c(n,s,M,\delta,r)\varepsilon^{-\frac{s}{n-s}}.
    \end{align}
\end{lemma}
The ``distribution of mass'' method dates back to the seminal work of
Chen and Li \cite{Chenli} and has been widely used in the study of improved
inequalities and prescribing curvature problems
(see, for instance, \cite{Dja,Ndi,Sun} and the references therein). To the best of our knowledge, previous work has mainly focused on the first-order case (i.e. $s=1$) and the Hilbert case (i.e. $n=2s$).
For the first-order case (see, for example, \cite{Chenli,Sun}), the proofs typically rely on a truncation argument originating in \cite{Aub2}.
For the Hilbert case (see, for example, \cite{Dja,Ndi}), the proofs are usually based on suitable eigenspace decompositions.
However, these methods do not extend to the general setting.

To establish Lemma \ref{lem6.1}, we draw on Fontana's proof of \eqref{MT3} in \cite{Fon}, whose ideas trace back to \cite{Ada}.
We use the Green representation of $u\in W^{s,\frac{n}{s}}(M)$ satisfying $\int_M u \de\mu =0$:
\begin{align}\label{Green}
    u(x)=\int_{M}G_s(x,y)\cdot\Delta^\frac{s}{2}u(y)\de\mu(y)=:J_s(\Delta^\frac{s}{2}u)(x).
\end{align}
Here $G_s$ is a scalar function when $s$ is even, and a vector field when $s$ is odd.
For any $v\in L^\frac{n}{s}(M)$, we define
\begin{align}\label{Riesz0}
    I_s(v)(x)=\int_{M}|G_s(x,y)||v(y)|\de\mu(y).
\end{align}
From \cite[Theorems 1.8 and 1.9]{Fon} and Young's inequality, we obtain the following Riesz-type Moser--Trudinger--Onofri inequality:
\begin{align}\label{Riesz}
     \log\int_M e^{nI_s(v)}\,\mathrm{d}\mu
    \leq \alpha_{s,n}
    \norm{v}_{L^{\frac{n}{s}}(M)}^{\frac{n}{s}}
    + C(M,s).
\end{align}
By \eqref{Green} and \eqref{Riesz0}, inequality \eqref{Riesz} implies \eqref{MTO0}.
Our main idea is to reduce the estimates for \eqref{MTO0} to those for \eqref{Riesz}, which can then be handled using classical cut-off techniques similar to those in \cite[Theorem 3.3]{Sun}.
\begin{proof}[Proof of Lemma \ref{lem6.1}]
    For simplicity, let $c_0$ denote a positive constant depending only on $n$, $s$, $M$, $\delta$, and $r$ and allowed to vary from line to line.

    Choose $\eta_i\in C^\infty(M)$ $(1\leq i\leq k)$ such that
    $0\leq \eta_i\leq 1$, $\eta_i|_{\widetilde{\Omega}_i}=1$, where
    $\widetilde{\Omega}_i:=\{x\in M:\;\text{dist}(x,\Omega_i)\leq \frac{r}{3}\}$,
    $|\nabla \eta_i|\leq c_0 r^{-1}$, and
    $\text{supp}(\eta_i)\cap \text{supp}(\eta_j)=\emptyset$ whenever $i\neq j$. For each $1\leq i\leq k$, equations \eqref{Green} and \eqref{distMTO1} give
    \begin{align}\label{dd1}
        \int_{M} e^{nu}\de\mu
        \leq&\; \delta^{-1} \int_{\Omega_i} e^{nu}\de\mu\nonumber\\
        =&\;\delta^{-1} \int_{\Omega_i} e^{nJ_s(\Delta^\frac{s}{2}u)}\,\mathrm{d}\mu\nonumber\\
        =&\;\delta^{-1} \int_{\Omega_i} e^{nJ_s(\eta_i\Delta^\frac{s}{2}u)
        +nJ_s((1-\eta_i)\Delta^\frac{s}{2}u)}\,\mathrm{d}\mu.
    \end{align}
    By \cite[Theorem 1.8]{Fon} and the properties of $\eta_i$, for any $x\in\Omega_i$,
    \begin{align}\label{dd2}
        J_s((1-\eta_i)\Delta^\frac{s}{2}u)(x)
        \leq c_0
        \norm{\Delta^\frac{s}{2}u}_{L^1(M)}.
    \end{align}
    Using \eqref{Riesz0}, \eqref{Riesz}, \eqref{dd1}, and \eqref{dd2}, we obtain
    \begin{align}\label{dd3}
        \int_{M} e^{nu}\de\mu
        \leq e^{C(M,s)}\delta^{-1}
        \exp\!\left(
        \alpha_{s,n}
        \norm{\eta_i\Delta^{\frac{s}{2}} u}_{L^{\frac{n}{s}}(M)}^{\frac{n}{s}}
        +c_0
        \norm{\Delta^\frac{s}{2}u}_{L^\frac{n}{s}(M)}
        \right).
    \end{align}
    By Young's inequality,
    \begin{align}\label{dd4}
        c_0
        \norm{\Delta^\frac{s}{2}u}_{L^\frac{n}{s}(M)}
        \leq
        \varepsilon
        \norm{\Delta^\frac{s}{2}u}_{L^\frac{n}{s}(M)}^\frac{n}{s}
        +c_0^{\frac{n}{n-s}}\varepsilon^{-\frac{s}{n-s}}.
    \end{align}
    Moreover, there exists an index $1\leq i\leq k$ such that
    \begin{align}\label{dd5}
        \norm{\eta_i\Delta^{\frac{s}{2}} u}_{L^{\frac{n}{s}}(M)}^{\frac{n}{s}}
        \leq
        \frac{1}{k}
        \norm{\Delta^{\frac{s}{2}} u}_{L^{\frac{n}{s}}(M)}^{\frac{n}{s}}.
    \end{align}
    Combining \eqref{dd3}, \eqref{dd4}, and \eqref{dd5} yields \eqref{distMTO2}.
\end{proof}
Combining Remark \ref{rmk4.4} with Lemma \ref{lem6.1} gives the following result.
\begin{theorem}\label{thm6.2}
Assume that $s,m,n\in\N_+$, $1\leq s<n$, and $0<\varepsilon<1$. Then, for any function $u\in W^{s,\frac{n}{s}}(\S^n)$ satisfying $
\int_{\S^n} u\,\mathrm{d}\mu = 0$ and
\begin{align}\label{or4}
    \left|\int_{\S^n} f e^{n u}\,\mathrm{d}\mu\right|
    \leq
    |T(f)|
    \quad
    \text{for all } f\in \PP_{m,n},
\end{align}
where $T$ is a linear functional on $\PP_{m,n}$, the following inequality holds:
\begin{align}\label{main4}
     \log \int_{\S^n} e^{n u}\,\mathrm{d}\mu
    \leq
    \left(
        \frac{\alpha_{s,n}}{\Theta(m,0,n)}
        +\varepsilon
    \right)
    \norm{\Delta^{\frac{s}{2}} u}_{L^{\frac{n}{s}}(\S^n)}^{\frac{n}{s}}
    +c_\varepsilon,
\end{align}
where the constant $c_\varepsilon$ can be chosen as
\begin{align}\label{cst5}
    c_\varepsilon
    =
    \log\big(\norm{T}_\infty+1\big)
    +C(m,n,s)\varepsilon^{-\frac{s}{n-s}},
\end{align}
and $\norm{T}_\infty$ is defined in \eqref{eq5.beta}.
\end{theorem}

\begin{proof}
Take $\varepsilon_0=\varepsilon_0(m,n)$ to be the positive constant in Remark \ref{rmk4.4}. If
\begin{align}\label{yy1}
    \int_{\S^n}e^{nu}\de\mu
    \leq
    \norm{T}_\infty\varepsilon_0^{-1},
\end{align}
then \eqref{main4} holds provided that
\[
c_\varepsilon
\geq
\log\big(\norm{T}_\infty+1\big)
-\log \varepsilon_0.
\]

Assume now that \eqref{yy1} does not hold. Define
\[
\de\sigma
=
\frac{e^{nu}\de\mu}
{\int_{\S^n}e^{nu}\,\mathrm{d}\mu}.
\]
Then $\sigma$ is a Borel probability measure on $\S^n$. Moreover, by \eqref{or4}, it satisfies the relaxed moment condition in Remark \ref{rmk4.4}. Therefore, we can find $N_*=\Theta(m,0,n)$ points $x_1,\cdots,x_{N_*}\in\S^n$ such that $d(x_i,x_j)\geq R=R(m,n)>0$ whenever $i\neq j$, and $\sigma(B(x_i,\frac{R}{200}))\geq \delta(m,n)>0$ for every $i$.

Applying Lemma \ref{lem6.1} to the domains $\Omega_i:=B(x_i,\frac{R}{200})$ $(1\leq i\leq N_*)$, we obtain \eqref{main4} from \eqref{distMTO2}. The estimate \eqref{cst5} follows from the preceding argument.
\end{proof}
We conclude this subsection by proving Corollary \ref{cor1.9}.
\begin{proof}[Proof of Corollary \ref{cor1.9}]
    For simplicity, let $C$ denote a constant depending only on $n,s,m$ and allowed to vary from line to line.

    From Theorem \ref{thm1.7}, for any $0<\varepsilon<1$,
    \begin{align}\label{vv1}
        \log \int_{\S^n} e^{n u}\,\mathrm{d}\mu
        \le
        \Big(\frac{\alpha_{s,n}}{\Theta(m,0,n)}+\varepsilon\Big)
        \norm{\Delta^{\frac{s}{2}} u}_{L^{\frac{n}{s}}(\S^n)}^{\frac{n}{s}}
        + C\varepsilon^{-\frac{s}{n-s}}.
    \end{align}
    It suffices to consider the case $\norm{\Delta^{\frac{s}{2}} u}_{L^{\frac{n}{s}}(\S^n)}>1$. Setting
    \begin{align}
        \varepsilon=\norm{\Delta^{\frac{s}{2}} u}_{L^{\frac{n}{s}}(\S^n)}^{1-\frac{n}{s}},\nonumber
    \end{align}
    we derive from \eqref{vv1} that
    \begin{align}\label{vv2}
        \log \int_{\S^n} e^{n u}\,\mathrm{d}\mu
        \le \frac{\alpha_{s,n}}{\Theta(m,0,n)} \norm{\Delta^{\frac{s}{2}} u}_{L^{\frac{n}{s}}(\S^n)}^{\frac{n}{s}}
        + C\norm{\Delta^{\frac{s}{2}} u}_{L^{\frac{n}{s}}(\S^n)}.
    \end{align}
    Choosing $q=\frac{n}{n-s}$,
    \begin{align}
        x=\left( \frac{\alpha_{s,n}}{\Theta(m,0,n)} \right)^{1-\frac{s}{n}} \norm{\Delta^{\frac{s}{2}} u}_{L^{\frac{n}{s}}(\S^n)}^{\frac{n}{s}-1}\nonumber
    \end{align}
    and
    \begin{align}
        y=C,\nonumber
    \end{align}
    we see from \eqref{bb3} that the right-hand side of \eqref{vv2} is bounded by $(x+y)^q$. Hence,
    \begin{align}
        \log^{\frac{1}{q}}\Big( \frac{1}{\omega_n}\int_{\S^n} e^{n u}\,\mathrm{d}\mu \Big) \le x+y,\nonumber
    \end{align}
    which gives \eqref{scaline2}.
\end{proof}
\subsection{Second-moment refinements}\label{sec6.2}
In this subsection, we revisit the recent work \cite{Chenzhang}. For $u\in W^{1,2}(\S^2)$, set $\Lambda(u)=(\Lambda_{ij}(u))_{i,j}$, where, for $1\leq i,j\leq 3$,
\begin{align}
    \Lambda_{ij}(u)=\frac{\int_{\S^2} e^{2 u(x)}\big(x_ix_j-\frac{1}{3}\delta_{ij}\big)\,\mathrm{d}\mu}{\int_{\S^2} e^{2 u}\,\mathrm{d}\mu}.\nonumber
\end{align}
Define $\norm{\Lambda(u)}^2=\sum_{i,j}\Lambda_{ij}(u)^2$, which measures the deviation of the second-order moments from zero. In \cite[Theorem 1.1]{Chenzhang}, Chen and Zhang proved that if $u$ satisfies $\int_{\S^2}u\,\mathrm{d}\mu=0$ and \eqref{or2} with $m=1$ and $n=2$, and if $\norm{\Lambda(u)}^2\leq \frac{2}{3}-\delta$ for some positive constant $\delta$, then, for any $0<\varepsilon<1$,
\begin{align}\label{mm1}
    \log \int_{\S^2} e^{2 u}\,\mathrm{d}\mu
    \leq
    \left(
        \frac{\alpha_{1,2}}{3}
        +\varepsilon
    \right)\norm*{\nabla u}_{L^2(\S^2)}^2+c_{\varepsilon,\delta}
\end{align}
for some positive constant $c_{\varepsilon,\delta}$. If, in addition, $\norm{\Lambda(u)}^2\leq \frac{1}{6}-\delta$, then
\begin{align}\label{mm2}
    \log \int_{\S^2} e^{2 u}\,\mathrm{d}\mu
    \leq
    \left(
        \frac{\alpha_{1,2}}{4}
        +\varepsilon
    \right)\norm*{\nabla u}_{L^2(\S^2)}^2+c_{\varepsilon,\delta}.
\end{align}

Note that $\Theta(1,0,2)=2$ and $\Theta(2,0,2)=4$ (see Section \ref{sec2}). Hence \eqref{mm1} represents an intermediate case between \eqref{main2} with $m=1,n=2$ and \eqref{main2} with $m=2,n=2$, while \eqref{mm2} improves \eqref{main2} by relaxing the condition \eqref{or2} in the case $m=2,n=2$. Moreover, both \eqref{mm1} and \eqref{mm2} are \textbf{almost sharp}, and the upper bounds $\frac{2}{3}-\delta$ and $\frac{1}{6}-\delta$ cannot be replaced by any constants larger than $\frac{2}{3}$ and $\frac{1}{6}$, respectively.

The proofs of \eqref{mm1} and \eqref{mm2} in \cite{Chenzhang} rely on a careful application of the concentration--compactness principles developed in \cite{Cha2,Hang}. Here, we give a new quantitative proof of these results and further generalize them to the higher-order inequality \eqref{MTO0} on $\S^n$ for $n\geq 2$.

We first introduce some notation. Let $\sigma$ be a Radon probability measure on $\S^n$. Define $\Lambda^{(n)}(\sigma)=(\Lambda^{(n)}_{ij}(\sigma))_{i,j}$, where for $1\leq i,j\leq n+1$,
\begin{align}
    \Lambda^{(n)}_{ij}(\sigma)=\int_{\S^n}\left(x_ix_j-\frac{1}{n+1}\delta_{ij}\right)\,\mathrm{d}\sigma.\nonumber
\end{align}
Define $\norm{\Lambda^{(n)}(\sigma)}^2=\sum_{i,j}\Lambda^{(n)}_{ij}(\sigma)^2$. For any function $u\in W^{s,\frac{n}{s}}(\S^n)$ with $\int_{\S^n}u\,\mathrm{d}\mu=0$, set $\Lambda^{(n)}(u)=\Lambda^{(n)}(\sigma_u)$, where
\begin{align}
    \mathrm{d}\sigma_u=\frac{e^{n u}\,\mathrm{d}\mu}{\int_{\S^n} e^{n u}\,\mathrm{d}\mu}.\nonumber
\end{align}
From Section \ref{sec2}, we know that $\Theta(1,0,n)=2$ and $\Theta(2,0,n)=n+2$. For any integer $2\leq k\leq n+2$, define
\begin{align}
    g_n(k)=\inf\left\{\norm{\Lambda^{(n)}(\nu)}^2:\nu\in\M_{1,n}^c,\ \norm{\nu}_0\leq k\right\}.\nonumber
\end{align}
Recall that $\norm{\nu}_0$ denotes the cardinality of the support of $\nu$. It was shown in \cite{Chenzhang} that $g_2(2)=\frac{2}{3}$, $g_2(3)=\frac{1}{6}$, and $g_2(4)=0$. The following lemma determines $g_n(k)$ in all dimensions.

\begin{lemma}\label{lem6.3}
For any integer $2\leq k\leq n+2$, we have
\[
    g_n(k)
    =
    \frac{1}{k-1}-\frac{1}{n+1}
    =
    \frac{n-k+2}{(k-1)(n+1)} .
\]
Moreover, equality is attained if and only if $\nu$ is the uniform measure
supported on the vertices of a regular $(k-1)$-simplex contained in a
$(k-1)$-dimensional linear subspace of $\mathbb R^{n+1}$. In particular, the
minimizing measure is unique up to orthogonal transformations of
$\mathbb R^{n+1}$.
\end{lemma}

\begin{proof}
Let
\[
    \nu=\sum_{a=1}^{s}w_a\delta_{x_a},
    \qquad 1\leq s\leq k,
\]
be an admissible measure, where $x_a\in\mathbb S^n$, $w_a>0$, and
$\sum\limits_{a=1}^sw_a=1$. Since $\nu\in\M^c_{1,n}$, we have
$\sum\limits_{a=1}^sw_ax_a=0$.

Define the second moment matrix
\[
    M=\int_{\mathbb S^n}x\otimes x\,\mathrm{d}\nu
    =
    \sum_{a=1}^{s}w_ax_ax_a^T .
\]
Then $M$ is symmetric and positive semidefinite. Moreover,
$\operatorname{tr}M=\sum\limits_{a=1}^sw_a|x_a|^2=1$. Since
\[
    \Lambda^{(n)}(\nu)=M-\frac1{n+1}I_{n+1},
\]
we obtain
\[
    \|\Lambda^{(n)}(\nu)\|^2
    =
    \operatorname{tr}(M^2)-\frac1{n+1}.
\]

It remains to estimate $\operatorname{tr}(M^2)$. Let
$r=\operatorname{rank}M$. Since
\[
\operatorname{Ran}(M)\subset
\operatorname{span}\{x_1,\ldots,x_s\},
\]
we have
\[
r\leq \dim\operatorname{span}\{x_1,\ldots,x_s\}.
\]
On the other hand, the relation $\sum\limits_{a=1}^{s}w_ax_a=0$
is a nontrivial linear dependence relation among the vectors
$x_1,\ldots,x_s$. Therefore,
\[
r\leq \dim\operatorname{span}\{x_1,\ldots,x_s\}\leq s-1\leq k-1.
\]

Let $\lambda_1,\ldots,\lambda_r$ be the nonzero eigenvalues of $M$. Since
$\sum\limits_{i=1}^{r}\lambda_i=\operatorname{tr}M=1$, the Cauchy--Schwarz
inequality gives
\[
\operatorname{tr}(M^2)
=
\sum_{i=1}^{r}\lambda_i^2
\geq
\frac1r
\left(\sum_{i=1}^{r}\lambda_i\right)^2
=
\frac1r .
\]
Combining this with $r\leq k-1$, we obtain
\[
\operatorname{tr}(M^2)\geq\frac1{k-1}.
\]
Consequently,
\[
\|\Lambda^{(n)}(\nu)\|^2
\geq
\frac1{k-1}-\frac1{n+1},
\]
which proves
\[
g_n(k)\geq
\frac1{k-1}-\frac1{n+1}.
\]

We next show that this lower bound is attained. Since $k\leq n+2$,
there exists a $(k-1)$-dimensional linear subspace
$V\subset\mathbb R^{n+1}$ and a regular $(k-1)$-simplex
$v_1,\ldots,v_k\in V\cap\mathbb S^n$ satisfying
\[
\sum_{a=1}^{k}v_a=0,
\qquad
\langle v_a,v_b\rangle=-\frac1{k-1}
\quad(a\neq b).
\]
Let $\nu_*=\frac1k\sum_{a=1}^{k}\delta_{v_a}$. Then $\nu_*\in\M^c_{1,n}$. Moreover, the vertices of a regular simplex
satisfy the identity
\[
\frac1k\sum_{a=1}^{k}v_av_a^T
=
\frac1{k-1}P_V ,
\]
where $P_V$ is the orthogonal projection onto $V$. Hence
\[
M_*=\frac1{k-1}P_V
\]
and therefore
\[
\operatorname{tr}(M_*^2)
=
\frac1{(k-1)^2}\operatorname{tr}(P_V)
=
\frac1{k-1}.
\]
Thus
\[
\|\Lambda^{(n)}(\nu_*)\|^2
=
\frac1{k-1}-\frac1{n+1},
\]
which proves the desired formula for $g_n(k)$.

It remains to characterize the equality case. Suppose that an admissible
measure $\nu$ attains the above lower bound. Then equality must hold in the
Cauchy--Schwarz inequality, and hence all nonzero eigenvalues of $M$ are
equal. Since
\[
\operatorname{tr}(M^2)=\frac1{k-1},
\]
we must have $\operatorname{rank}M=k-1$.
In particular, the support of $\nu$ contains exactly $k$ points.

Let
\[
A=(\sqrt{w_1}x_1,\ldots,\sqrt{w_k}x_k)
\]
and set
\[
q=(\sqrt{w_1},\ldots,\sqrt{w_k})^T .
\]
The barycenter condition implies $Aq=0$. Since
$\operatorname{rank}A=k-1$, we have $\ker A=\operatorname{span}\{q\}$. Moreover, $AA^T=M$
has exactly $k-1$ nonzero eigenvalues, all equal to $1/(k-1)$. Hence
$A^TA$ has the same nonzero eigenvalues and kernel
$\operatorname{span}\{q\}$, which yields
\[
A^TA=\frac1{k-1}(I_k-qq^T).
\]

Comparing the diagonal entries, we obtain
\[
w_i
=
\frac1{k-1}(1-w_i),
\]
and hence $w_i=\frac1k$
for every $i$. Comparing the off-diagonal entries gives
\[
\sqrt{w_iw_j}\langle x_i,x_j\rangle
=
-\frac1{k-1}\sqrt{w_iw_j},
\]
and therefore
\[
\langle x_i,x_j\rangle=-\frac1{k-1},
\qquad i\neq j .
\]
Thus the support points form a regular $(k-1)$-simplex, and the measure is
the uniform measure on its vertices. Conversely, we have already shown that
every such measure attains the lower bound. The proof is complete.
\end{proof}

We are now ready to state our main result.
\begin{theorem}
Assume that $s,n,k\in\N_+$, $1\leq s<n$, $3\leq k\leq n+2$, $0<\delta<\frac{1}{(k-1)(k-2)}$, and $0<\varepsilon<1$. Then, for any function $u\in W^{s,\frac{n}{s}}(\S^n)$ satisfying $\int_{\S^n} u\,\mathrm{d}\mu = 0$, $\norm{\Lambda^{(n)}(u)}^2\leq g_n(k-1)-\delta$, and \eqref{or2} with $m=1$, we have
\begin{align}\label{main6}
     \log \int_{\S^n} e^{n u}\,\mathrm{d}\mu
    \leq
    \left(
        \frac{\alpha_{s,n}}{k}
        +\varepsilon
    \right)
    \norm{\Delta^{\frac{s}{2}} u}_{L^{\frac{n}{s}}(\S^n)}^{\frac{n}{s}}
    +C(n,s,\delta,k)\varepsilon^{-\frac{s}{n-s}}.
\end{align}
\end{theorem}
Standard constructions as in \cite{Cha2,Chenzhang,Hang,Sun} show that \eqref{main6} is \textbf{almost sharp} and that the upper bound $g_n(k-1)-\delta$ cannot be replaced by any constant larger than $g_n(k-1)$. More precisely, choose probability measures $\nu_k$ and $\nu_k'$ attaining $g_n(k)$ and $g_n(k-1)$, respectively. Following \cite[Propositions 2.3 and 2.4]{Chenzhang} and \cite[Theorem 3.3]{Hang}, one constructs test functions concentrating near the supports of $\nu_k$ and $\nu_k'$, respectively, and performs direct expansions. Since these computations are standard, we omit the details.

As in the proof of Theorem \ref{thm6.2}, it suffices to combine Lemma \ref{lem6.1} with the following result to deduce \eqref{main6}.

\begin{lemma}
Assume that $n,k\in\N_+$, $3\leq k\leq n+2$, and $0<\delta<\frac{1}{(k-1)(k-2)}$. If $\sigma$ is a Radon probability measure on $\S^n$ whose first moments vanish and such that $\norm{\Lambda^{(n)}(\sigma)}^2\leq g_n(k-1)-\delta$, then there exist positive constants $\gamma_1=\gamma_1(n)$ and $\gamma_2=\gamma_2(n)$ and $k$ open subsets $\Omega_1,\cdots,\Omega_k\subset \S^n$ satisfying $\mathrm{dist}(\Omega_i,\Omega_j)\geq \gamma_1\delta$ whenever $i\neq j$ and $\sigma(\Omega_i)\geq \gamma_2\delta^{n+1}$ for every $i$.
\end{lemma}
\begin{proof}
For simplicity, let $C_n$ denote a positive constant depending only on $n$ and allowed to vary from line to line. By Lemma~\ref{lem6.3}, we have $g_n(k-1)=\frac1{k-2}-\frac1{n+1}$. Set $M(\sigma)=\int_{\S^n}x\otimes x\,\mathrm{d}\sigma$. The inequality $\norm{\Lambda^{(n)}(\sigma)}^2\leq g_n(k-1)-\delta$ gives
\begin{align}\label{op0}
    \operatorname{tr}(M(\sigma)^2)
    \leq
    \frac1{k-2}-\delta .
\end{align}

We call a spherical ball $B(x,\gamma_1\delta)$ heavy if $\sigma(B(x,\gamma_1\delta))\geq\gamma_2\delta^{n+1}$. Here $\gamma_1=\gamma_1(n)$ and $\gamma_2=\gamma_2(n)$ are positive constants to be determined. Choose a maximal family of heavy balls $B(x_1,\gamma_1\delta),\ldots,B(x_\ell,\gamma_1\delta)$ whose centers satisfy $|x_i-x_j|\geq4\gamma_1\delta$ whenever $i\neq j$.

If $\ell\geq k$, then any $k$ of these balls give the desired conclusion.
Suppose, for contradiction, that $\ell\leq k-1$.

Define
\[
    U=\bigcup_{i=1}^{\ell}B(x_i,6\gamma_1\delta).
\]
Choose a maximal $\gamma_1\delta$-separated set
$\{y_1,\ldots,y_N\}$ in $\S^n\setminus U$. A standard covering estimate gives $N\leq C_n(\gamma_1\delta)^{-n}$. Moreover, the balls $B(y_j,\gamma_1\delta)$ are not heavy, since otherwise they could
be added to the maximal family. Hence $\sigma(B(y_j,\gamma_1\delta))<\gamma_2\delta^{n+1}$. Because these balls cover $\S^n\setminus U$, we obtain
\begin{align}\label{op1}
    \sigma(\S^n\setminus U)
    \leq C_n(\gamma_1\delta)^{-n}\gamma_2\delta^{n+1}
    =C_n\gamma_1^{-n}\gamma_2\delta.
\end{align}
Let $E_1,\ldots,E_\ell$ be a disjoint Borel partition of $U$ satisfying $E_i\subset B(x_i,6\gamma_1\delta)$
and set
\[
    p_i=\sigma(E_i),\qquad
    P=\sum_{i=1}^{\ell}p_i .
\]
Then, by \eqref{op1}, we have $P\geq 1-C_n\gamma_1^{-n}\gamma_2\delta$.

Define
\[
    a=\sum_{i=1}^{\ell}p_ix_i,
    \qquad
    C=\sum_{i=1}^{\ell}p_ix_ix_i^T .
\]
Since the first moments of $\sigma$ vanish, we can write
\[
    a
    =
    \sum_i\int_{E_i}(x_i-x)\,\mathrm{d}\sigma
    -
    \int_{\S^n\setminus U}x\,\mathrm{d}\sigma .
\]
Therefore, $|a|\leq 6\gamma_1\delta+C_n\gamma_1^{-n}\gamma_2\delta$.

Similarly, we have
\[
    M(\sigma)-C
    =
    \sum_i\int_{E_i}
    (xx^T-x_ix_i^T)\,\mathrm{d}\sigma
    +
    \int_{\S^n\setminus U}xx^T\,\mathrm{d}\sigma .
\]
By H\"older's inequality we have $\|M(\sigma)-C\|
    \leq C_n\gamma_1\delta+C_n\gamma_1^{-n}\gamma_2\delta$.

Now let $A$ be the matrix whose columns are $\sqrt{p_i}x_i$, and let
\[
    v=(\sqrt{p_1},\ldots,\sqrt{p_\ell})^T .
\]
Then
\[
    AA^T=C,\qquad Av=a,\qquad |v|^2=P .
\]
Set $A^TA=G$ and let $\lambda_1$ be the smallest eigenvalue. By the Rayleigh quotient and the estimates of $|a|$ and $P$ above, we obtain
\begin{align}\label{op2}
    \lambda_1
    \leq
    \frac{v^TGv}{|v|^2}
    =
    \frac{|a|^2}{P}\leq e:= \frac{(6\gamma_1\delta+C_n\gamma_1^{-n}\gamma_2\delta)^2}{1-C_n\gamma_1^{-n}\gamma_2\delta}.
\end{align}

Since $\ell\leq k-1$, after removing the smallest eigenvalue, $G$ has at
most $k-2$ nonzero eigenvalues. By \eqref{op2}, their sum is at least $P-e$. Therefore,
by the Cauchy--Schwarz inequality, we have
\[
    \operatorname{tr}(C^2)
    =
    \operatorname{tr}(G^2)
    \geq
    \frac{(P-e)^2}{k-2}.
\]
Since both $M(\sigma)$ and $C$ are positive semidefinite and
have trace bounded by one,
\[
    |\operatorname{tr}(M(\sigma)^2)-\operatorname{tr}(C^2)|
    \leq
    2\|M(\sigma)-C\|
    \leq
    2C_n\gamma_1\delta+2C_n\gamma_1^{-n}\gamma_2\delta.
\]
Combining the preceding two estimates yields
\begin{align}
    \operatorname{tr}(M(\sigma)^2)\geq \frac{(P-e)^2}{k-2}-(2C_n\gamma_1\delta+2C_n\gamma_1^{-n}\gamma_2\delta).\nonumber
\end{align}
Using the estimate $P\geq 1-C_n\gamma_1^{-n}\gamma_2\delta$ and the expression for $e$ in \eqref{op2}, we deduce, provided that both $\gamma_1$ and $\gamma_1^{-n}\gamma_2$ are sufficiently small, that
\begin{align}
     \operatorname{tr}(M(\sigma)^2)\geq \frac{1}{k-2}-\frac{\delta}{2},\nonumber
\end{align}
which contradicts the estimate \eqref{op0}. The proof is thus complete.
\end{proof}

An argument similar to the proof of Corollary \ref{cor1.9} also gives the following \textbf{sharp} improved inequality at a different scale.

\begin{corollary}
Assume that $s,n,k\in\N_+$, $1\leq s<n$, $3\leq k\leq n+2$, and $0<\delta<\frac{1}{(k-1)(k-2)}$. Then, for any function $u\in W^{s,\frac{n}{s}}(\S^n)$ satisfying $\int_{\S^n} u\,\mathrm{d}\mu = 0$, $\norm{\Lambda^{(n)}(u)}^2\leq g_n(k-1)-\delta$, and \eqref{or2} with $m=1$, we have
\begin{align}
     \log^{1-\frac{s}{n}}\Bigg( \frac{1}{\omega_{n}}\int_{\S^n} e^{n u}\,\mathrm{d}\mu\Bigg)
        \leq
        \left(
            \frac{\alpha_{s,n}}{k}
        \right)^{1-\frac{s}{n}}
        \norm{\Delta^{\frac{s}{2}} u}_{L^{\frac{n}{s}}(\S^n)}^{\frac{n}{s}-1}
    +C(n,s,\delta,k).\nonumber
\end{align}
\end{corollary}

\section{Improved Sobolev-type inequalities under symmetry constraints}\label{sec7}
In this section, we study how symmetry can improve Sobolev-type embeddings. The first systematic investigation of this subject on $\R^n$ was carried out by Lions in the celebrated work \cite{Lions0}. Here we focus on smooth closed connected manifolds. In general, two types of improvement are possible for functions invariant under the action of a group $G$.

\vskip0.1in

First, if every orbit of $G$ is infinite, one may obtain Sobolev-type embeddings into higher $L^p$ spaces and, in particular, compact embeddings. Second, if $G$ has finite orbits, embeddings into higher $L^p$ spaces are no longer expected; one may instead seek improved embedding constants. We focus on the second type. We refer to \cite[Chapter 5]{Heb0}, \cite[Chapter 9]{Heb1}, and \cite{Fag,Heb3} for further discussion, including detailed treatments of noncompact manifolds.

\vskip0.1in

We begin with some notation and preliminaries; see also \cite[Section 9.1]{Heb1}. Let $(M,g)$ be a smooth closed manifold of dimension $n\geq 2$, and let $\operatorname{Iso}(M)$ denote its isometry group. Suppose that $G$ is a compact subgroup of $\operatorname{Iso}(M)$. For $x\in M$, define its orbit by
\begin{align}
    O_G^x=\{gx:g\in G\}.\nonumber
\end{align}
We set
\[
k_G=\min_{x\in M}|O_G^x|,
\]
where $|O_G^x|$ denotes the cardinality of $O_G^x$. Throughout this section, we assume $k_G<+\infty$. For $s\in\N_+$ and $1\leq p\leq+\infty$, define
\begin{align}
    W_G^{s,p}(M)=\{u\in W^{s,p}(M):u\circ g=u \text{ for all } g\in G\}.\nonumber
\end{align}

In the pioneering work \cite[Theorem 5.8]{Heb0}, Hebey and Vaugon established the following improved \textbf{almost sharp} Sobolev inequality. If $1\leq p<n$ and $u\in W_G^{1,p}(M)$, then, for every $0<\varepsilon<1$, there exists a positive constant $b_\varepsilon$ such that
\begin{align}\label{main7}
    \norm*{u}_{L^{p^*}(M)}^p
    \leq
    \Big(
        \frac{S_{n,p}^p}{k_G^{\frac{p}{n}}}
        +\varepsilon
    \Big)
    \norm*{\nabla u}_{L^p(M)}^p
    +b_\varepsilon\norm*{u}_{L^{p}(M)}^p.
\end{align}
In the special case $p=2$, they also proved the \textbf{sharp} inequality (see \cite[Theorem 5.10]{Heb0})
\begin{align}\label{main8}
    \norm*{u}_{L^{p^*}(M)}^2
    \leq
    \frac{S_{n,2}^2}{k_G^{\frac{2}{n}}}
    \norm*{\nabla u}_{L^2(M)}^2
    +b\norm*{u}_{L^{2}(M)}^2
\end{align}
for some positive constant $b$. Their proof relies on a careful construction of a partition of unity and a detailed local analysis along the $G$-orbits.

\vskip0.1in

In \cite[Theorem 2]{Fag}, Faget established a first-order Moser--Trudinger--Onofri analogue of \eqref{main7}. More precisely, if $u\in W_G^{1,n}(M)$ and $\int_{M}u\,\mathrm{d}\mu=0$, then, for every $0<\varepsilon<1$, there exists a positive constant $c_\varepsilon$ such that
\begin{align}\label{main9}
    \log \int_{M} e^{n u}\,\mathrm{d}\mu
    \leq
    \left(
        \frac{\alpha_{1,n}}{k_G}
        +\varepsilon
    \right)
    \norm{\nabla u}_{L^n(M)}^{n}+c_\varepsilon.
\end{align}
This inequality is \textbf{almost sharp}. Its proof in \cite{Fag} combines the techniques developed in \cite{Heb0} with a truncation method introduced in \cite{Aub2}.

\vskip0.1in

The strategies used to prove Theorems \ref{thm5.1} and \ref{thm6.2} yield the following generalizations of \eqref{main8} and \eqref{main9}.
\begin{theorem}\label{thm7.1}
Assume $1\leq p<n$ and $u\in W_G^{1,p}(M)$. Suppose that the \textbf{sharp} Sobolev inequality \eqref{Sob3} holds. If $1\leq p\leq 2$, then
\begin{align}
    \norm*{u}_{L^{p^*}(M)}^p
    \leq
        \frac{S_{n,p}^p}{k_G^{\frac{p}{n}}}
    \norm*{\nabla u}_{L^p(M)}^p
    +b_G(M,p)\norm*{u}_{L^{p}(M)}^p,\nonumber
\end{align}
whereas if $2<p<n$, then, for every $0<\varepsilon<1$,
\begin{align}
    \norm*{u}_{L^{p^*}(M)}^p
    \leq
        \Big(
        \frac{S_{n,p}^p}{k_G^{\frac{p}{n}}}
        +\varepsilon
    \Big)
    \norm*{\nabla u}_{L^p(M)}^p
    +b_G(M,p)\varepsilon^{1-\frac{p}{2}}\norm*{u}_{L^{p}(M)}^p.\nonumber
\end{align}
\end{theorem}

\vskip0.1in

\begin{theorem}\label{thm7.2}
Assume $s\in\N_+$, $1\leq s<n$, and $u\in W_G^{s,\frac{n}{s}}(M)$ with $\int_M u\,\mathrm{d}\mu=0$. Then, for every $0<\varepsilon<1$,
\begin{align}
    \log \int_{M} e^{n u}\,\mathrm{d}\mu
    \leq
    \left(
        \frac{\alpha_{s,n}}{k_G}
        +\varepsilon
    \right)
    \norm{\Delta^{\frac{s}{2}} u}_{L^{\frac{n}{s}}(M)}^{\frac{n}{s}}
    +c_G(M,s)\varepsilon^{-\frac{s}{n-s}}.\nonumber
\end{align}
\end{theorem}
Theorem \ref{thm7.1} gives a general criterion for the validity of improved \textbf{sharp} Sobolev inequalities under symmetry constraints. In particular, \eqref{main8} is a special case of Theorem \ref{thm7.1}, because Hebey--Vaugon \cite{Heb2} proved that, when $p=2$, the \textbf{sharp} Sobolev inequality \eqref{Sob3} holds on every smooth closed manifold. To our knowledge, Theorem \ref{thm7.1} is the first general sharpness result since the work of Hebey--Vaugon \cite{Heb2}. Theorem \ref{thm7.2}, in turn, extends \eqref{main9} to the higher-order setting.

\vskip0.1in

The key ingredients in the proofs are the following three lemmas.
\begin{lemma}\label{lem7.3}
Let $\sigma$ be a Borel probability measure on $M$. If $\sigma$ is invariant under the action of $G$, then there exists $r_G>0$ such that
\[
\sigma(B(x,r_G))\leq\frac{1}{k_G}
\]
for any $x\in M$.
\end{lemma}

\begin{proof}
Fix $x\in M$. Since $|O_G^x|\geq k_G$, choose distinct points $y_1=x,y_2,\dots,y_{k_G}$ in $O_G^x$. There exists $r_x>0$ such that
\[
B(y_i,r_x)\cap B(y_j,r_x)=\emptyset \quad \text{for } i\neq j.
\]
Because $\sigma$ is invariant under the action of $G$, we have
\[
\sigma(B(y_i,r_x))=\sigma(B(y_j,r_x))
\]
for all $i\neq j$. Hence
\[
\sigma(B(x,r_x))\leq \frac{1}{k_G}.
\]
Covering $M$ by finitely many balls $B(x_i,\frac{r_{x_i}}{2})$, $1\leq i\leq I$, and setting $r_G:=\min\limits_i\{\frac{r_{x_i}}{2}\}$ gives the desired conclusion.
\end{proof}
\begin{lemma}\label{lem7.4}
Assume $N>0$, $1\leq p<n$, $r>0$, and $\sigma$ is a Borel probability measure on $M$ satisfying
\begin{align}\label{est1}
    \sigma\big(B(x,r)\big)
        \leq
        \frac{1}{N}
        \quad \text{for all } x\in M.
\end{align}
Then there exist finitely many points $x_j\in M$, $1\leq j\leq J$, and an $\ell^p$-partition of unity $\{\eta_j\}_{j=1}^{J}$
such that $\operatorname{supp}\eta_j
\subset B(x_j,r)$,
\begin{align}
    \sum_{j=1}^{J}\eta_j^p=1,\qquad \sum_{j=1}^{J}|\nabla\eta_j|^p\leq Cr^{-p},\qquad \sum_{j=1}^{J}\beta_j^{\frac{n-p}{n}}
\geq N^{\frac{p}{n}},
\end{align}
where $\beta_j:=
\int_{M}\eta_j^{p^*}\,d\sigma$ and $C=C(M,p)$.
\end{lemma}
\begin{proof}
Throughout the proof, $C$ denotes a positive constant depending only on $M$ and $p$, whose value may vary from line to line. Take a maximal $\frac{r}{8}$-separated set
$\{x_j\}_{j=1}^{J}$ of
$M$. We can find a family of smooth functions $\{\varphi_j\}_{j=1}^J$ such that $\supp \varphi_j\subset B(x_j,r)$, $0\leq \varphi_j\leq 1$, $\varphi_j=1$ on $B(x_j,\frac{r}{2})$, and $|\nabla \varphi_j|\leq Cr^{-1}$. Define
\[
\eta_j
=
\frac{\varphi_j}
{\left(\sum\limits_{k=1}^J\varphi_k^p\right)^{1/p}} .
\]
Then $\sum\limits_{j=1}^J\eta_j^p=1$. By the choice of $\{x_j\}_{j=1}^{J}$, every point $y\in M$ belongs to at most $C$ of the balls $B(x_i,r)$. Therefore,
$\sum\limits^J_{j=1}|\nabla\eta_j|^p\leq Cr^{-p}$.

\vskip0.1in

It remains to prove the estimate involving $\beta_j$. By H\"older's inequality and \eqref{est1},
\begin{align}
    \int_{M}\eta_j^p\,\de\sigma\leq\left(
\int_{M}\eta_j^{p^*}\,\de\sigma
\right)^{p/p^*}
\sigma(\operatorname{supp}\eta_j)^{1-p/p^*}\leq \beta_j^{p/p^*}N^{-p/n}.\nonumber
\end{align}
Summing over $j$ and using $\sum\limits_{j=1}^J\eta_j^p=1$, we obtain
\begin{align}
    \sum_{j=1}^{J}\beta_j^{\frac{n-p}{n}}\geq N^{\frac{p}{n}},\nonumber
\end{align}
which completes the proof.
\end{proof}
\begin{lemma}\label{lem7.5}
    Let $\sigma$ be a Borel probability measure on $M$. If $\sigma$ is invariant under the action of $G$, then there exist open subsets $\Omega_1,\dots,\Omega_{k_G}\subset M$ and positive constants $R=R(M,G)$ and $\delta=\delta(M,G)$ such that $\operatorname{dist}(\Omega_i,\Omega_j)\geq R$ whenever $i\neq j$, and $\sigma(\Omega_i)\geq \delta$ for every $1\leq i\leq k_G$.
\end{lemma}
\begin{proof}
    By Lemma \ref{lem7.3}, the argument is similar to that used in the proof of Lemma \ref{lem4.3}; we therefore omit the details.
\end{proof}
Replacing Proposition \ref{prop5.7} and Lemma \ref{lem4.3} by Lemmas \ref{lem7.4} and \ref{lem7.5}, respectively, in the proofs of Theorems \ref{thm5.1} and \ref{thm6.2} yields Theorems \ref{thm7.1} and \ref{thm7.2}. We omit the details.

\vskip0.1in

Arguments similar to those used in the proofs of Corollaries \ref{cor1.8} and \ref{cor1.9} also give the following \textbf{sharp} improved inequalities at different scales.
\begin{corollary}
Assume that $2<p<n$ and $u\in W_G^{1,p}(M)$. If \eqref{Sob3} holds, then
\begin{align}
    \norm*{u}_{L^{p^*}(M)}^2
    \leq
    \Big(
        \frac{S_{n,p}}{k_G^{\frac{1}{n}}}
    \Big)^{2}
    \norm*{\nabla u}_{L^p(M)}^{2}
    +D(M,p,G)\norm*{u}_{L^{p}(M)}^{2}.\nonumber
\end{align}
\end{corollary}

\begin{corollary}
Let $s\in\N_+$ with $1\leq s<n$, and let $u\in W_G^{s,\frac{n}{s}}(M)$ satisfy $\int_M u\,\mathrm{d}\mu=0$. Then
\begin{align}
    \log^{1-\frac{s}{n}}\Bigg( \frac{1}{V(M)}\int_{M} e^{n u}\,\mathrm{d}\mu\Bigg)
    \leq
    \left(
        \frac{\alpha_{s,n}}{k_G}
    \right)^{1-\frac{s}{n}}
    \norm{\Delta^{\frac{s}{2}} u}_{L^{\frac{n}{s}}(M)}^{\frac{n}{s}-1}
    +D(M,s,G),
\end{align}
where $V(M)$ denotes the volume of $(M,g)$.
\end{corollary}
We conclude this section with a general criterion for the failure of \textbf{sharp} improved Sobolev inequalities under symmetry constraints. To our knowledge, it gives the first non-sharp case of an improved Sobolev inequality under symmetry constraints.
\begin{theorem}\label{thm7.8}
    Assume $2<p<\sqrt{n}$ and that there exists a point $x\in M$ such that $k_G=|O_G^x|$ and the scalar curvature at $x$ is positive. Then there is no finite constant $b$ such that
    \begin{align}
    \norm*{u}_{L^{p^*}(M)}^p
    \leq
        \frac{S_{n,p}^p}{k_G^{\frac{p}{n}}}\norm*{\nabla u}_{L^p(M)}^p
    +b\norm*{u}_{L^{p}(M)}^p\nonumber
\end{align}
holds for every $u\in W_G^{1,p}(M)$.
\end{theorem}
\begin{proof}
    The proof is analogous to that of Theorem \ref{thm1.2}. The proof of \cite[Theorem 1]{Dru} gives a sequence of nonnegative functions $u_k\in C_c^\infty(M)$ such that $\supp u_k\subset B(x,k^{-1})$, each $u_k(y)$ depends only on $d_g(x,y)$, and
    \begin{align}
        \frac{
\|u_k\|_{L^{p^*}(M)}^{p}
-
S_{n,p}^{p}\|\nabla u_k\|_{L^p(M)}^{p}
}
{\|u_k\|_{L^p(M)}^{p}}
\rightarrow+\infty.\nonumber
    \end{align}
    Write $O_G^x=\{g_1x,\dots,g_{k_G}x\}$, where $g_i\in G$, and set $v_{k,i}=u_k\circ (g_i)^{-1}$. Then $\supp v_{k,i}\subset B(g_ix,k^{-1})$, and, for all sufficiently large $k$, $\supp v_{k,i}\cap\supp v_{k,j}=\emptyset$ whenever $i\neq j$. For such $k$, define $w_k=\sum\limits^{k_G}_{i=1}v_{k,i}\in W_G^{1,p}(M)$. Then
    \begin{align}
        \frac{
\|w_k\|_{L^{p^*}(M)}^{p}
-
S_{n,p}^{p}k_G^{-\frac{p}{n}}\|\nabla w_k\|_{L^p(M)}^{p}
}
{\|w_k\|_{L^p(M)}^{p}}
=\frac{
\|u_k\|_{L^{p^*}(M)}^{p}
-
S_{n,p}^{p}\|\nabla u_k\|_{L^p(M)}^{p}
}
{k_G^{\frac{p}{n}}\|u_k\|_{L^p(M)}^{p}}
\rightarrow+\infty,\nonumber
    \end{align}
    which completes the proof.
\end{proof}

\vskip0.2in
\noindent{\Large\textbf{Declarations}}
\vskip0.2in

\medskip
\noindent\textbf{Declaration of generative AI and AI-assisted technologies.}
During the preparation of this work, the authors used ChatGPT (OpenAI)
to assist with language editing and locating relevant references, and
to assist in searching for explicit examples for the case
$m=1$, $5\le n\le 8$, in Theorem~4.6.
The suggested references and examples, together with the associated
mathematical calculations, were independently checked and verified
by the authors. The authors reviewed and edited the content as needed
and take full responsibility for the contents of the article.
 
\vskip0.12in
\smallskip
\noindent\textbf{Conflict of interest.}\quad On behalf of all authors, the corresponding author states that there is no conflict of interest.

\vskip0.12in
\smallskip
\noindent\textbf{Data availability statement.}\quad No data were used for the research described in this article.


\addcontentsline{toc}{section}{References}\begin{thebibliography}{99}

\bibitem{An} An, J., Dou, J., Han, Y.: Almost sharp weighted Sobolev trace inequalities in the unit ball under constraints. Commun. Contemp. Math. 26, no. 6, Paper No. 2350025, 29 pp (2024)

\bibitem{Ada} Adams, D. R.: A sharp inequality of J. Moser for higher order derivatives. Ann. of Math. (2) 128, no. 2, 385--398 (1988)

\bibitem{Adi} Adimurthi, Druet, O.:
Blow-up analysis in dimension 2 and a sharp form of Trudinger--Moser inequality. Comm. Partial Differential Equations 29, nos. 1--2, 295--322 (2004)

\bibitem{Aub1} Aubin, T.: Problemes isop\'erim\'etriques et espaces de Sobolev. J. Differential Geom. 11, 573--598 (1976)

\bibitem{Aub2} Aubin, T.: Meilleures constantes dans le th\'eoreme d'inclusion de Sobolev et un th\'eoreme de Fredholm non lin\'eaire pour la transformation conforme de la courbure scalaire. J. Functional Analysis 32, no. 2, 148--174 (1979)

\bibitem{Aub3} Aubin, T.: Some nonlinear problems in Riemannian geometry. Springer Monographs in Mathematics. Springer-Verlag, Berlin, xviii+395 pp (1998)

\bibitem{Bec} Beckner, W.: Sharp Sobolev inequalities on the sphere and the Moser--Trudinger inequality. Ann. of Math. (2) 138, no. 1, 213--242 (1993)

\bibitem{Bia} Bianchi, G., Egnell, H.: A note on the Sobolev inequality. J. Funct. Anal. 100, 18--24 (1991)

\bibitem{Bra} Branson, T. P., Chang, S.-Y. A., Yang, P. C.:
Estimates and extremals for zeta function determinants on four-manifolds. Comm. Math. Phys. 149, no. 2, 241--262 (1992)



\bibitem{Cha2} Chang, S.-Y. A., Hang, F.: Improved Moser--Trudinger--Onofri inequality under constraints. Comm. Pure Appl. Math. 75, no. 1, 197--220 (2022)

\bibitem{Cha3} Chang, S.-Y. A., Yang, P. C.: Prescribing Gaussian curvature on $\S^2$. Acta Math. 159, no. 3--4, 215--259 (1987)

\bibitem{Cha4} Chang, S.-Y. A., Yang, P. C.: A perturbation result in prescribing scalar curvature on $\S^n$.
Duke Math. J. 64, no. 1, 27--69 (1991)


\bibitem{Chenli} Chen, W., Li, C.: Prescribing Gaussian curvatures on surfaces with conical singularities. J. Geom. Anal. 1, no. 4, 359--372 (1991)

\bibitem{Chenwei} Chen, X., Wei, W., Wu, N.: Almost sharp Sobolev trace inequalities in the unit ball under constraints. Adv. Math. 459, Paper No. 110023, 88 pp (2024)

\bibitem{Chenzhang} Chen, X., Zhang, S.: Constrained Moser--Trudinger--Onofri inequality and a uniqueness criterion for the mean field equation. Calc. Var. Partial Differential Equations 64, no. 9, Paper No. 296, 30 pp (2025)

\bibitem{de} de la Harpe, P., Pache, C.: Cubature formulas, geometrical designs, reproducing kernels, and Markov operators. Infinite groups: geometric, combinatorial and dynamical aspects, 219--267,
Progr. Math., 248, Birkhäuser, Basel (2005)

\bibitem{Del} Delsarte, P., Goethals, J. M., Seidel, J. J.: Spherical codes and designs. Geometriae Dedicata 6, no. 3, 363--388 (1977)

\bibitem{Dja} Djadli, Z., Malchiodi, A.: Existence of conformal metrics with constant Q-curvature. Ann. of Math. (2) 168, no. 3, 813--858 (2008)

\bibitem{Dru} Druet, O.: Optimal Sobolev inequalities of arbitrary order on compact Riemannian manifolds. J. Funct. Anal. 159, no. 1, 217--242 (1998)

\bibitem{Fag} Faget, Z.: Best constants in Sobolev inequalities on Riemannian manifolds in the presence of symmetries. Potential Anal. 17, no. 2, 105--124 (2002)

\bibitem{Fon} Fontana, L.:
Sharp borderline Sobolev inequalities on compact Riemannian manifolds.
Comment. Math. Helv. 68, no. 3, 415--454 (1993)

\bibitem{Fon2} Fontana, L., Morpurgo, C.: Adams inequalities on measure spaces, Adv. Math. 226, no. 6, 5066--5119 (2011)

\bibitem{FJWX}
Frank, R.L., Jin, T., Wang, H.-Y., Xiong, J.:
Sharp Sobolev inequalities under first moment constraints.
arXiv preprint (2026), arXiv:2610.20001.



\bibitem{Guihu1} Gui, C., Hu, Y., Xie, W.: Improved Beckner's inequality for axially symmetric functions on $\S^n$.
J. Funct. Anal. 282, no. 5, Paper No. 109335, 47 pp (2022)

\bibitem{Guihu2} Gui, C., Hu, Y., Xie, W.: Improved Beckner's inequality for axially symmetric functions on $\S^4$. Rev. Mat. Iberoam. 40, no. 1, 355--388 (2024)

\bibitem{Guili} Gui, C., Li, T., Wei, J., Ye, Z.: Sharp Beckner's inequalities for axially symmetric functions on $\S^6$ and $\S^8$. Adv. Math. 480, part A, Paper No. 110487, 64 pp (2025)

\bibitem{Guili1} Gui, C., Li, T., Wei, J., Ye, Z.: A positive answer to the generalized Chang-Yang conjecture on $\S^N$. arXiv preprint (2026), arXiv:2608.13497

\bibitem{Guili2} Gui, C., Li, T., Wei, J., Ye, Z.: Sharp Beckner's Inequalities for Axially Symmetric Functions on $\S^N$. arXiv preprint (2026), arXiv:2608.11126

\bibitem{Guimo} Gui, C., Moradifam, A.: The sphere covering inequality and its applications. Invent. Math. 214, no. 3, 1169--1204 (2018)

\bibitem{Hang} Hang, F.: A remark on the concentration compactness principle in critical dimension. Comm. Pure Appl. Math. 75, no. 10, 2245--2278 (2022)

\bibitem{Hangwang} Hang, F., Wang, X.: Improved Sobolev inequality under constraints. Int. Math. Res. Not., no. 14, 10822--10857 (2022)


\bibitem{Heb0} Hebey, E.: Sobolev spaces on Riemannian manifolds. Lecture Notes in Mathematics, 1635. Springer-Verlag, Berlin. x+116 pp (1996)

\bibitem{Heb1} Hebey, E.: Nonlinear analysis on manifolds: Sobolev spaces and inequalities. Courant Lecture Notes in Mathematics, 5. New York University, Courant Institute of Mathematical Sciences, New York; American Mathematical Society, Providence, RI, x+309 pp (1999)

\bibitem{Heb2} Hebey, E., Vaugon, M.: Meilleures constantes dans le th\'eoreme d'inclusion de Sobolev. Ann. Inst. H. Poincar\'e C Anal. Non Lin\'eaire 13, no. 1, 57--93 (1996)

\bibitem{Heb3} Hebey, E., Vaugon, M.: Sobolev spaces in the presence of symmetries.
J. Math. Pures Appl. (9) 76, no. 10, 859--881 (1997)

\bibitem{Ho} Ho, P. T.: Improved $L^p$-Folland--Stein--Sobolev inequality under constraints. Arch. Math. (Basel) 120, no. 5, 521--532 (2023)

\bibitem{Hor1} H\"ormander, L.: The Analysis of Linear Partial Differential Operators I: Distribution
Theory and Fourier Analysis, 2nd ed., Springer, Berlin, 1990.

\bibitem{Hor2} H\"ormander, L.: The Analysis of Linear Partial Differential Operators III: Pseudo-Differential Operators, Springer, Berlin, 1985; reprint, 2007.

\bibitem{Hu} Hu, Y., Dou, J.:
Improved Hardy--Littlewood--Sobolev inequality on $\S^n$ under constraints. Acta Math. Sin. (Engl. Ser.) 39, no. 11, 2149--2163 (2023)

\bibitem{Li} Li, T., Wei, J., Ye, Z.: On sharp Beckner's inequality for axially symmetric functions on $\S^4$. Math. Res. Lett. 31, no. 5, 1523--1550 (2024)

\bibitem{Lions0} Lions, P.-L.: Sym\'etrie et compacit\'e dans les espaces de Sobolev. J. Functional Analysis 49, no. 3, 315--334 (1982)

\bibitem{Lions1} Lions, P.-L.:
The concentration-compactness principle in the calculus of variations. The limit case. I.
Rev. Mat. Iberoamericana 1, no. 1, 145--201 (1985)

\bibitem{Lions2} Lions, P.-L.:
The concentration-compactness principle in the calculus of variations. The limit case. II.
Rev. Mat. Iberoamericana 1, no. 2, 45--121 (1985)

\bibitem{Moser} Moser, J.: A sharp form of an inequality by N. Trudinger. Indiana Univ. Math. J. 20, 1077--1092 (1970)

\bibitem{Ndi} Ndiaye, C. B.: Constant $Q$-curvature metrics in arbitrary dimension. J. Funct. Anal. 251, no. 1, 1--58 (2007)

\bibitem{Onofri} Onofri, E.: On the positivity of the effective action in a theory of random surfaces. Comm. Math. Phys. 86, no. 3, 321--326 (1982)

\bibitem{Put} Putterman, E.: Cubature formulas and Sobolev inequalities. arXiv preprint (2020), arXiv:2012.08109

\bibitem{See} Seeley, R. T.: Complex powers of an elliptic operator, in Singular Integrals, Proc. Sympos. Pure Math., vol. 10, American
Mathematical Society, Providence, RI, 1967, pp. 288--307.

\bibitem{Shu} Shubin, M. A.: Pseudodifferential Operators and Spectral Theory, 2nd ed., Springer,
Berlin, 2001.

\bibitem{Stein} Stein, E. M., Weiss, G.: Introduction to Fourier analysis on Euclidean spaces. Princeton Mathematical Series, No. 32. Princeton University Press, Princeton, NJ. x+297 pp (1971)

\bibitem{Sun} Sun, J., Zou, W.: Some improved almost sharp inequalities under constraints. Int. Math. Res. Not., no. 16, Paper No. rnaf242, 27 pp (2025)

\bibitem{Sze} Szego, G.: Orthogonal Polynomials, 4th ed., American Mathematical Society, Colloquium Publications, Vol. 23, (1975)

\bibitem{Tal} Talenti, G.: Best constant in Sobolev inequality.
Ann. Mat. Pura Appl. (4) 110, 353--372 (1976)

\bibitem{Tin} Tintarev, K.: Trudinger--Moser inequality with remainder terms. J. Funct. Anal. 266, no. 1, 55--66 (2014)

\bibitem{Tru} Trudinger, N. S.:
On imbeddings into Orlicz spaces and some applications.
J. Math. Mech. 17, 473--483 (1967)

\bibitem{Wang} Wang, G., Ye, D.: A Hardy--Moser--Trudinger inequality. Adv. Math. 230, no. 1, 294--320 (2012)

\bibitem{WYY-Dolbeault}
Wei, J., Yan, Z., Ye, Z.:
On Dolbeault's conjecture.
arXiv preprint (2026), arXiv:2610.14359.

\bibitem{WYY-Sobolev}
Wei, J., Yan, Z., Ye, Z.:
Sharp improved Sobolev inequalities on $\S^n$.
arXiv preprint (2026), arXiv:2610.23069.

\bibitem{Yan1} Yan, Z.: Improved higher-order Sobolev inequalities on CR sphere. J. Funct. Anal. 284, no. 10, Paper No. 109890, 34 pp (2023)

\bibitem{Yan2} Yan, Z.: Improved Sobolev inequalities on CR sphere. Calc. Var. Partial Differential Equations 64, no. 2, Paper No. 70, 21 pp (2025)

\bibitem{Yang} Yang, Y.: On a sharp inequality of L. Fontana for compact Riemannian manifolds. Manuscripta Math. 157, no. 1--2, 51--79 (2018)


\end{thebibliography}
\end{document}